\documentclass[12pt]{article}%
\usepackage{amsfonts}
\usepackage{amsmath}
\usepackage{amssymb}
\usepackage{graphicx}
\usepackage{color}
\usepackage{setspace}%
\providecommand{\U}[1]{\protect\rule{.1in}{.1in}}
\newtheorem{theorem}{Theorem}

\newtheorem{corollary}[theorem]{Corollary}

\newtheorem{example}[theorem]{Example}

\newtheorem{lemma}[theorem]{Lemma}

\newtheorem{proposition}[theorem]{Proposition}
\newtheorem{remark}[theorem]{Remark}

\newenvironment{proof}[1][Proof]{\noindent\textbf{#1.} }{\ \rule{0.5em}{0.5em}}
\begin{document}

\title{Rank Reversals in Efficient Vectors for Pairwise Comparisons}
\author{Susana Furtado\thanks{Corresponding author. Email: sbf@fep.up.pt.
orcid.org/0000-0003-0395-5972. The work of this author is funded by FCT -
Funda\c{c}\~{a}o para a Ci\^{e}ncia e a Tecnologia, I.P., through national
funds, under the project UID/04561/2025,
https://doi.org/10.54499/UID/04561/2025.}\\CEMS.UL and Faculdade de Economia e Gest\~{a}o\\Universidade do Porto\\Rua Dr. Roberto Frias\\4200-464 Porto, Portugal
\and Charles R. Johnson \thanks{Email: crjmatrix@gmail.com. }\\225 West Tazewells Way\\Williamsburg, VA 23185, United States}
\maketitle

\begin{abstract}
Studied are the orders among alternatives given by efficient vectors for a
reciprocal matrix $A$. The focus is upon: 1) existence of efficient vectors
with a given order and 2) discrepancies between the ranking for a pair of
alternatives resulting from an efficient vector and from the corresponding
entry of $A$, called rank reversals. We use matrix intervals in which the
consistent matrices built from efficient vectors lie. The existence of
efficient vectors (or those from constituent convex subsets) with a rank
reversal at a given position is characterized. The minimum number of rank
reversals for a reciprocal matrix among all weight vectors is described and it
is shown that there is an efficient vector that has this number of rank
reversals. A way to construct such a vector is given. The order transitive and
the uniformly ordered reciprocal matrices are shown to be, respectively, the
matrices for which there is an efficient vector with no rank reversals and the
matrices for which all efficient vectors have no rank reversals. The effects
on rank reversals of addition and deletion of alternatives are also discussed.
Using numerical experiments, the frequencies of rank reversals in certain
formulaic weight vectors are given.

\end{abstract}



\section{Introduction}

A \emph{pairwise comparison}, or \emph{reciprocal}, matrix is an $n$-by-$n$
entrywise positive matrix $A=[a_{ij}]$ in which $a_{ji}=\frac{1}{a_{ij}},$
$1\leq i,j\leq n$. The diagonal entries are all $1$, and the off-diagonal
entries represent pairwise ratio comparisons among $n$ alternatives. Let
$\mathcal{PC}_{n}$ denote the set of all such matrices. If $A\in
\mathcal{PC}_{n}$ satisfies, in addition, $a_{ij}a_{jk}=a_{ik}$, $1\leq
i,j,k\leq n$, then $A$ is further said to be \emph{consistent}.

In a variety of ranking models, a weight vector $w=\left[
\begin{array}
[c]{ccc}%
w_{1} & \cdots & w_{n}%
\end{array}
\right]  ^{T}\in\mathbb{R}_{+}^{n}$ (the set of positive $n$-vectors) is to be
deduced from $A\in\mathcal{PC}_{n}$, so that $\frac{w_{i}}{w_{j}}\approx
a_{ij}$. If $A$ is consistent, then $A=ww^{(-T)}=\left[  \frac{w_{i}}{w_{j}%
}\right]  $, for some $w\in\mathbb{R}_{+}^{n}$. Here, $w^{(-T)}$ is the
entrywise inverse of the transpose of $w$. In this case, $w$ is (projectively)
the natural weight vector. However, consistency of comparisons seldom occurs
in practice. For this reason, it is important to obtain a weight vector from
an inconsistent reciprocal matrix. We denote the set of such matrices in
$\mathcal{PC}_{n}$ by $\mathcal{PC}_{n}^{0}$. Several methods for obtaining a
weight vector have been proposed and discussed, see surveys
\cite{Choo,Golany1993}.

We ask that the consistent matrix $W=ww^{(-T)}$ approximates $A\in
\mathcal{PC}_{n}$ in such a way that no other $v\in\mathbb{R}_{+}^{n}$ better
approximates $A$ via $V=vv^{(-T)}$, i.e. $\left\vert A-V\right\vert
\leq\left\vert A-W\right\vert ,$ entrywise, implies that $v$ is proportional
to $w.$ Such a vector $w$ is called \emph{efficient} (or \emph{Pareto optimal}
from an economics perspective). We denote the set of all efficient vectors for
$A$ by $\mathcal{E}(A)$, which will be projectively infinite if $A\in
\mathcal{PC}_{n}^{0}$. For a consistent matrix, $\mathcal{E}(A)$ is just the
positive multiples of any of its columns, and thus is projectively a single
vector. Two nice methods to describe all efficient vectors for each
$A\in\mathcal{PC}_{n}^{0}$ have recently been developed; $\mathcal{E}(A)$ as a
union of a finite number of convex sets \cite{FJ3}, and an inductice
construction of $\mathcal{E}(A)$ \cite{FJ4}. The former is described in
Section \ref{sunion} and used throughout, and the latter is described and used
in Section \ref{sadddelete}.

We will say that $w\in\mathbb{R}_{+}^{n}$ has an $i,j$ \emph{(strict) rank
reversal} for $A\in\mathcal{PC}_{n}^{0}$ if (1) $w_{i}>w_{j}$ when $a_{ij}<1$,
or (2) $w_{i}<w_{j}$ when $a_{ij}>1$. This means that comparison of
alternatives $i$ and $j$ coming from $w$ gives a different order from the one
coming from $a_{ij}$. It is widely agreed that the number of rank reversals is
an important consideration when choosing a weight vector. Some authors have
studied rank reversals in a reciprocal matrix relative to some often used
formulaic weight vectors, such as the right Perron eigenvector
\cite{saaty1977} and the geometric mean of the columns \cite{William}. The
occurrence of rank reversals when alternatives are added or deleted, as well
as for completions of partial reciprocal matrices, have been discussed. See
for example
\cite{Bana,Barzilai-Cs,BEL,Csato2017,CasatoLP,Csato,Fara2,Maleki-Cs,MAZ,Saaty-Cs,TuWu,Wang1-Cs,Wang2-Cs,wang}
and references therein for an overview of the results available in the literature.

Given $A\in\mathcal{PC}_{n}^{0}$, finding a vector with few rank reversals for
$A$ among efficient vectors is desirable. Though there has been discussion
about rank reversals for weight vectors, little is known for efficient weight
vectors. The latter is our purpose. Recent developments in the theory of
efficient vectors make this feasible.

According to our definition, if $a_{ij}=1$ or $w_{i}=w_{j}$, there is no $i,j$
rank reversal for $A$ at $w$. For simplicity, we will confine our attention to
matrices $A\in\mathcal{PC}_{n}^{0}$ with no off-diagonal entries equal to $1$
and often to vectors $w\in\mathbb{R}_{+}^{n}$ with distinct entries. We denote
the set of such matrices $A$ by $\mathcal{PC}_{n}^{(1)}$ and call such a
vector $w$ a \emph{strict} vector. Restricted to such matrices and vectors,
our definition coincides with the one in \cite{Golany1993}. The order of the
entries in a strict vector is called a \emph{strict order}.

There is a natural, one-way, relationship between matrices in $\mathcal{PC}%
_{n}^{(1)}$ and "tournaments". A \emph{tournament} on vertices $1,\ldots,n$ is
a directed graph with, for each pair $i\neq j$, exactly one of the edges
$(i,j)$ or $(j,i)$ \cite{Moon}. The adjacency matrix $T$ of a tournament is an
$n$-by-$n$ $0,1$ matrix such that $T+T^{T}=J-I$. Here, $J$ is the matrix of
ones and $I$ the identity matrix. The tournament matrix $T(A)$ associated with
reciprocal matrix $A=[a_{ij}]\in\mathcal{PC}_{n}^{(1)}$ has a $1$ in the $i,j$
entry if and only if $a_{ij}>1$. Thus, it identifies the pattern of entries
$>1$ or $<1$ in $A$. There is a long history of structural investigation of
tournaments, including the ranking of contestants. Since this, at times,
resembles the problem of weight vectors for reciprocal matrices, we may
exploit tournament theory by applying it to $T(A)$. For example, the
definition of rank reversal for $A$ at a positive vector can be stated in
terms of $T(A)$ (Section \ref{sbest}). Note that a tournament is just a
combinatorial version of the quantitative notion of a reciprocal matrix $A$,
and there is no way to reconstruct the richer notion of $A$ from $T(A)$.

A set $P_{n}\subseteq\mathbb{R}_{+}^{n}$ is called $i,j$ \emph{order
unambiguous}, for $1\leq i,j\leq n$, $i\neq j$, if $\frac{p_{i}}{p_{j}}\geq1$
for all $p\in P_{n}$, or $\frac{p_{i}}{p_{j}}\leq1$ for all $p\in P_{n}$.
Otherwise, it is said to be $i,j$ \emph{order ambiguous}. This concept is
closely related to the nonexistence of a rank reversal in (a subset of)
$\mathcal{E}(A)$, as we will see.

In the discussion here about rank reversals in efficient vectors, we collect
known facts from scattered places, give new results based on recent work, and
develop some new ideas. The structure of the paper is as follows.

In Section \ref{sbackground} we survey useful results from recent work which
made some of the advances here possible. For $A\in\mathcal{PC}_{n}^{0}$, we
present a decomposition of $\mathcal{E}(A)$ into finitely many convex sets,
each associated with a certain full cycle in $A$. Also, matrix intervals, in
which the consistent matrices based upon vectors in these sets lie, are
presented. In Section \ref{Sorder} the existence of an efficient vector for a
reciprocal matrix with a given order is investigated. In the $3$-by-$3$ case,
necessary and sufficient conditions for the existence of such a vector are
given. In Section \ref{s3} we give necessary and sufficient conditions for the
existence of an $i,j$ order ambiguity in each of the constituent convex sets
of $\mathcal{E}(A)$, as well as in $\mathcal{E}(A)$. In Section \ref{srrev},
we relate $i,j$ order ambiguity in $\mathcal{E}(A)$, or in a constituent
convex set, with the existence of a vector in the set with an $i,j$ rank
reversal for $A$. Conditions for the nonexistence of an $i,j$ rank reversal
among all vectors in a constituent convex set of $\mathcal{E}(A)$ are
provided. It is shown that no vector in that set has a rank reversal in all
positions along the associated cycle. The existence of a vector with no rank
reversals along the cycle is characterized. In Section \ref{srankmin} we
describe the weight vectors that give the fewest rank reversals for
$A\in\mathcal{PC}_{n}^{(1)}$, and show that there is such a vector that is
efficient. An explicit way of obtaining this efficient vector is given. In
Section \ref{stransitive} we make some observations concerning the reciprocal
matrices for which there is a weight vector with no rank reversals. We call
these matrices \emph{order transitive}. In Section \ref{sUO} we focus on a
special class of order transitive reciprocal matrices, called \emph{uniformly
ordered}. These are the matrices for which all the efficient vectors have the
same order or, equivalently, no efficient vector has a rank reversal. We also
consider the case in which all the efficient vectors in a convex constituent
subset of $\mathcal{E}(A)$, not necessarily in all of $\mathcal{E}(A)$, have
the same order. In both cases, we give several different characterizations of
these occurrences. In Section \ref{sadddelete} we study the preservation of
cardinal rankings in efficient vectors when compared alternatives are added or
deleted. In Section \ref{snum} various formulaic weight vectors are studied
numerically relative to rank reversals for a reciprocal matrix. These vectors
include the often used right Perron eigenvector, the geometric mean of the
right and entrywise inverse left Perron eigenvector, and weighted geometric
means of the columns. The latter two have recently been shown to be
universally efficient \cite{FJEJOR,FJ2}, while the former may not be efficient
\cite{Boz,FJ5}. A summary of our major conclusions is given in Section
\ref{concl}.

\section{Background\label{sbackground}}

\subsection{Basic facts about efficiency}

A \emph{monomial matrix} is a product of a permutation and a positive diagonal matrix.

\begin{lemma}
\label{lsim}\cite{FJ1} If $A\in\mathcal{PC}_{n}$ and $S$ is an $n$-by-$n$
monomial matrix, then $SAS^{-1}\in\mathcal{PC}_{n}$ and is consistent if $A$
is. Moreover, $\mathcal{E}(SAS^{-1})=S\mathcal{E}(A)$.
\end{lemma}

A (directed) cycle $\alpha$ in a directed graph is a path (no repeated
vertices) with an additional edge from the last to the first vertex:
$i_{1}\rightarrow i_{2}\rightarrow i_{3}\rightarrow\cdots\rightarrow
i_{k}\rightarrow i_{k+1}=i_{1}$ (or $i_{1}i_{2}i_{3}\cdots i_{k}i_{1}$).
A\ path is called \emph{Hamiltonian} if it contains all vertices. If a cycle
contains a Hamiltonian path, it is called a \emph{Hamiltonian cycle} (an H-cycle).

\bigskip

For $A\in\mathcal{PC}_{n}$ and $w\in\mathbb{R}_{+}^{n}$, there is a graph
theoretic method to decide whether $w\in\mathcal{E}(A).$ Let $G(A,w)$ be the
digraph with vertex set $\{1,\ldots,n\}$ and a directed edge $(i,j)$ if and
only if $w_{i}\geq a_{ij}w_{j}$, $1\leq i\neq j\leq n$.

\begin{theorem}
\cite{blanq2006,FJ2,FJ3}\label{blanq} Suppose that $A\in\mathcal{PC}_{n}$ and
$w\in\mathbb{R}_{+}^{n}$. The following are equivalent:

\begin{enumerate}
\item[(i)] $w\in\mathcal{E}(A);$

\item[(ii)] $G(A,w)$ is a strongly connected digraph;

\item[(iii)] $G(A,w)$ contains a Hamiltonian cycle.
\end{enumerate}
\end{theorem}

\begin{remark}
\label{remspannintree}Suppose that $A\in\mathcal{PC}_{n}$. A labelled tree
(connected acyclic graph) $T$ on $n$ vertices may be viewed as a spanning tree
of the complete undirected graph with $n$ vertices. There is a unique
consistent matrix $ww^{(-T)}$ whose entries in positions given by the edges of
$T$ coincide with those of $A$ in the same positions (entries $a_{ij}$ and
$a_{ji}$ for the undirected edge $\{i,j\}$). By Theorem \ref{blanq}, vector
$w$ is efficient for $A$ because, if $\{i,j\}$ is an edge in $T$, then both
$(i,j)$ and $(j,i)$ are edges in $G(A,w),$ implying that $G(A,w)$ is strongly
connected. So, there is an element of $\mathcal{E}(A)$ associated with any
spanning tree of the complete graph. Examples are the columns of $A$
\cite{FJ1}, which correspond to the spanning trees that are labelled stars.
Another example is a vector $w$ such that $\frac{w_{\gamma_{i}}}%
{w_{\gamma_{i+1}}}=a_{\gamma_{i},\gamma_{i+1}}$ for all $i\in\{1,\ldots
,n-1\},$ in which $\gamma_{1}\cdots\gamma_{n}$ is a path on $n$ vertices in
the complete graph.
\end{remark}

\subsection{$\mathcal{E}(A)$ as a union of convex sets\label{sunion}}

For $A\in\mathcal{PC}_{n}^{0}$, the set $\mathcal{E}(A)$ is known to be
connected \cite{blanq2006,FJ4}. However, it is usually not convex
\cite{FJEJOR}. Next we describe $\mathcal{E}(A)$ as a union of a finite number
of constituent convex sets.

We identify the entries of $A=[a_{ij}]$ with the edges of the complete
directed graph on $n$ vertices, and view these interchangeably. Associated
with the cycle $\alpha=i_{1}i_{2}\cdots i_{k}i_{1}$ is the product of the
corresponding entries of $A$: $\alpha(A)=a_{i_{1}i_{2}}\cdots a_{i_{k}i_{1}}$.
If the cycle is Hamiltonian, its product in $A$ is called an \emph{H-cycle
product}. If we reverse the direction of the edges in a cycle $\alpha$, we
refer to the alternate cycle as $\operatorname*{rev}\alpha$. Of course
$\operatorname*{rev}\alpha(A)=\frac{1}{\alpha(A)},$ because $A$ is reciprocal.
So, at most one of the two can have product $<1$.

We denote by $\Gamma(A)$ the collection of all H-cycles with product $<1$ in
$A$.\ We have $\#\Gamma(A)\leq\frac{(n-1)!}{2},$ and $\#\Gamma(A)\geq1$ since
$A$ is inconsistent. If $A\in\mathcal{PC}_{3}^{0}$, then $\#\Gamma(A)=1$.

The cycles in $\Gamma(A)$ are important in determining $\mathcal{E}(A)$ for
$A\in\mathcal{PC}_{n}^{0}$. If $\gamma\in\Gamma(A)$, with $\gamma=\gamma
_{1}\cdots\gamma_{n}\gamma_{1}$, we denote by $\mathcal{E}_{\gamma}(A)$ the
set of vectors $w=\left[
\begin{array}
[c]{ccc}%
w_{1} & \cdots & w_{n}%
\end{array}
\right]  ^{T}\in\mathbb{R}_{+}^{n}$ such that
\begin{equation}
w_{\gamma_{1}}\geq a_{\gamma_{1}\gamma_{2}}w_{\gamma_{2}}\geq a_{\gamma
_{1}\gamma_{2}}a_{\gamma_{2}\gamma_{3}}w_{\gamma_{3}}\geq\cdots\geq\left(
{\displaystyle\prod\limits_{i=1}^{n-1}}
a_{\gamma_{i}\gamma_{i+1}}\right)  w_{\gamma_{n}}\geq\gamma(A)w_{\gamma_{1}}.
\label{in}%
\end{equation}
Notice that $\mathcal{E}_{\gamma}(A)$ just depends on the entries of $A$ along
$\gamma$.

\bigskip

Recently, the following description of $\mathcal{E}(A)$ was given.

\begin{theorem}
\cite{FJ3}\label{thunion} If $A\in\mathcal{PC}_{n}^{0}$,
\[
\mathcal{E}(A)=%
{\displaystyle\bigcup\limits_{\gamma\in\Gamma_{A}}}
\mathcal{E}_{\gamma}(A).
\]

\end{theorem}

Since the solution set of the linear inequalities (\ref{in}) in $\mathbb{R}%
_{+}^{n}$ is convex, this gives $\mathcal{E}(A)$ as the union of a finite
number of convex sets. Each set has infinitely many (projectively distinct) members.

\bigskip

We illustrate Theorem \ref{thunion} with an example that will serve multiple
purposes throughout the paper.

\begin{example}
\label{exbound0}Let%
\[
A=\left[
\begin{array}
[c]{cccc}%
1 & 2 & 2 & 16\\
\frac{1}{2} & 1 & 2 & 2\\
\frac{1}{2} & \frac{1}{2} & 1 & 2\\
\frac{1}{16} & \frac{1}{2} & \frac{1}{2} & 1
\end{array}
\right]  \in\mathcal{PC}_{4}^{(1)}.
\]
We have $\Gamma(A)=\{\alpha,\beta\},$ with $\alpha=12341$, $\beta=13241,$
$\alpha(A)=\frac{1}{2}$ and $\beta(A)=\frac{1}{8}.$ Note that $\mu(A)=1$ for
$\mu=12431$. By Theorem \ref{thunion}, $\mathcal{E}(A)=\mathcal{E}_{\alpha
}(A)\cup\mathcal{E}_{\beta}(A)$, with%
\[
\mathcal{E}_{\alpha}(A)=\left\{  w\in\mathbb{R}_{+}^{4}:w_{1}\geq2w_{2}%
\geq4w_{3}\geq8w_{4}\geq\frac{1}{2}w_{1}\right\}
\]
and%
\[
\mathcal{E}_{\beta}(A)=\left\{  w\in\mathbb{R}_{+}^{4}:w_{1}\geq2w_{3}\geq
w_{2}\geq2w_{4}\geq\frac{1}{8}w_{1}\right\}  .
\]

\end{example}

See \cite[Example 15]{FJ8} for a matrix $A\in\mathcal{PC}_{4}^{0}$ with
$\#\Gamma(A)=3$ (the generic case).

\subsection{Matrix intervals for consistent matrices coming from
$\mathcal{E}(A)$}

Let $A\in\mathcal{PC}_{n}^{0}$. Within $\gamma\in\Gamma_{A}$, for each pair
$i,j$, $i\neq j$, there is a directed path from $i$ to $j$. Denote the product
of the entries in $A$ corresponding to the edges in this path by $P_{A,\gamma
}(i,j)$. If $(i,j)$ is an edge in $\gamma$, then $P_{A,\gamma}(i,j)=a_{ij}$.
Also,%
\begin{equation}
P_{A,\gamma}(i,j)P_{A,\gamma}(j,i)=\gamma(A)<1. \label{Pijji}%
\end{equation}
We denote by $P_{A,\gamma}$ the $n$-by-$n$ matrix with $1$'s on the diagonal
and the $i,j$ entry $P_{A,\gamma}(i,j)$, $i\neq j$.

For a positive matrix $Q$, we denote by $Q^{(-T)}$ the transpose of its
entrywise inverse. For $w\in\mathbb{R}_{+}^{n}$, let $W=ww^{(-T)}$, the
consistent matrix built from $w$. Inequalities between matrices are to be
interpreted entrywise throughout.

\bigskip

Recently, the following was shown.

\begin{theorem}
\cite{FJ8}\label{Thintgamma} If $A\in\mathcal{PC}_{n}^{0}$, $\gamma\in
\Gamma(A)$ and $w\in\mathcal{E}_{\gamma}(A)$, then
\begin{equation}
P_{A,\gamma}\leq W\leq P_{A,\gamma}^{(-T)}. \label{int}%
\end{equation}
Moreover, any consistent matrix $W$ satisfying (\ref{int}) comes from some
$w\in\mathcal{E}_{\gamma}(A)$.
\end{theorem}

\bigskip Note that, because of (\ref{Pijji}), the interval in (\ref{int})
corresponding to any $i,j$ position, $i\neq j$, is not a singleton. Moreover,
$P_{A,\gamma}^{(-T)}=\frac{1}{\gamma(A)}P_{A,\gamma}+(1-\frac{1}{\gamma(A)}%
)I$. Since $W^{(-T)}=W,$ the second inequality in (\ref{int}) is equivalent to
the first. Also, each of these inequalities is equivalent to the set of
inequalities in (\ref{int}) restricted to the positions below (or above) the diagonal.

\bigskip

The collection of intervals (\ref{int}), with $\gamma\in\Gamma(A)$, for a
consistent matrix obtained from an efficient vector is a refinement of earlier
work \cite{FJ80} giving a single matrix interval that was strong, but not
quite as precise. This result can now be viewed as a corollary of Theorem
\ref{Thintgamma}.

Let $L_{A}$ denote the entrywise minimum (over $\gamma\in\Gamma(A)$) of
$P_{A,\gamma}$. Corresponding to (\ref{Pijji}), we also have
\begin{equation}
L_{A}(i,j)L_{A}(j,i)<1. \label{inLA}%
\end{equation}

\begin{corollary}
\cite{FJ80}\label{Thintgea} If $A\in\mathcal{PC}_{n}^{0}$ and $w\in
\mathcal{E}(A)$, then%
\begin{equation}
L_{A}\leq W\leq L_{A}^{(-T)}\text{.} \label{bLA}%
\end{equation}
Moreover, $L_{A}\leq A\leq L_{A}^{(-T)}$.
\end{corollary}

However, not all consistent $W$ satisfying (\ref{bLA}) are built from
efficient vectors.

\bigskip

Because of the next result, the matrix intervals given in Theorem
\ref{Thintgamma} and Corollary \ref{Thintgea} are the best possible.

\begin{theorem}
\cite{FJ80,FJ8}\label{rattain} If $A\in\mathcal{PC}_{n}^{0}$, $\gamma\in
\Gamma(A)$ and $1\leq i,j\leq n$, $i\neq j$, there are $w,v\in\mathcal{E}%
_{\gamma}(A)$ such that $\frac{w_{i}}{w_{j}}=P_{A,\gamma}(i,j)$ and
$\frac{v_{i}}{v_{j}}=\frac{1}{P_{A,\gamma}(j,i)}$ (the $i,j$ entry of
$P_{A,\gamma}^{(-T)}$). In particular, there are $w,v\in\mathcal{E}(A)$ such
that $\frac{w_{i}}{w_{j}}=L_{A}(i,j)$ and $\frac{v_{i}}{v_{j}}=\frac{1}%
{L_{A}(j,i)}$.
\end{theorem}

\bigskip Each of the vectors $w$ and $v$ in the first claim of the theorem can
be taken to be an extreme of $\mathcal{E}_{\gamma}(A)$.

\bigskip

Since both $\mathcal{E}(A)$ and $\mathcal{E}_{\gamma}(A)$ are connected, there
is a continuous path in each set between any two vectors in the set. Thus,
from Theorem \ref{rattain} it follows that, for any $x_{ij}\in\left[
L_{A}(i,j),\frac{1}{L_{A}(j,i)}\right]  $ ($x_{ij}\in\left[  P_{A,\gamma
}(i,j),\frac{1}{P_{A,\gamma}(j,i)}\right]  $) there is a $v\in\mathcal{E}(A)$
($v\in\mathcal{E}_{\gamma}(A)$) such that $\frac{v_{i}}{v_{j}}=x_{ij}$.

\begin{example}
\label{exbound02}Continue Example \ref{exbound0}. By Theorem \ref{Thintgamma},
we have $w\in\mathcal{E}_{\alpha}(A)$ if and only if
\begin{equation}
P_{A,\alpha}=\left[
\begin{array}
[c]{cccc}%
1 & 2 & 4 & 8\\
\frac{1}{4} & 1 & 2 & 4\\
\frac{1}{8} & \frac{1}{4} & 1 & 2\\
\frac{1}{16} & \frac{1}{8} & \frac{1}{4} & 1
\end{array}
\right]  \leq W\leq\left[
\begin{array}
[c]{cccc}%
1 & 4 & 8 & 16\\
\frac{1}{2} & 1 & 4 & 8\\
\frac{1}{4} & \frac{1}{2} & 1 & 4\\
\frac{1}{8} & \frac{1}{4} & \frac{1}{2} & 1
\end{array}
\right]  =P_{A,\alpha}^{(-T)}, \label{in1}%
\end{equation}
and $w\in\mathcal{E}_{\beta}(A)$ if and only if
\begin{equation}
P_{A,\beta}=\left[
\begin{array}
[c]{cccc}%
1 & 1 & 2 & 2\\
\frac{1}{8} & 1 & \frac{1}{4} & 2\\
\frac{1}{16} & \frac{1}{2} & 1 & 1\\
\frac{1}{16} & \frac{1}{16} & \frac{1}{8} & 1
\end{array}
\right]  \leq W\leq\left[
\begin{array}
[c]{cccc}%
1 & 8 & 16 & 16\\
1 & 1 & 2 & 16\\
\frac{1}{2} & 4 & 1 & 8\\
\frac{1}{2} & \frac{1}{2} & 1 & 1
\end{array}
\right]  =P_{A,\beta}^{(-T)}. \label{in2}%
\end{equation}
By Theorem \ref{thunion}, $w\in\mathcal{E}(A)$ if and only if at least one of
(\ref{in1}) or (\ref{in2}) holds. This implies (Corollary \ref{Thintgea})
that, if $w\in\mathcal{E}(A)$, then%
\[
L_{A}=\left[
\begin{array}
[c]{cccc}%
1 & 1 & 2 & 2\\
\frac{1}{8} & 1 & \frac{1}{4} & 2\\
\frac{1}{16} & \frac{1}{4} & 1 & 1\\
\frac{1}{16} & \frac{1}{16} & \frac{1}{8} & 1
\end{array}
\right]  \leq W\leq\left[
\begin{array}
[c]{cccc}%
1 & 8 & 16 & 16\\
1 & 1 & 4 & 16\\
\frac{1}{2} & 4 & 1 & 8\\
\frac{1}{2} & \frac{1}{2} & 1 & 1
\end{array}
\right]  =L_{A}^{(-T)}.
\]
Notice that the sets $\mathcal{E}_{\gamma}(A)$, $\gamma\in\Gamma(A)$, given in
Example \ref{exbound0}, can be written using entries of $P_{A,\gamma}$.\ For
example,
\[
\mathcal{E}_{\beta}(A)=\left\{  w\in\mathbb{R}_{+}^{4}:w_{1}\geq P_{A,\beta
}(1,3)w_{3}\geq P_{A,\beta}(1,2)w_{2}\geq P_{A,\beta}(1,4)w_{4}\geq
\beta(A)w_{1}\right\}  .
\]

\end{example}

\bigskip

Not only does Theorem \ref{Thintgamma} give another presentation of
$\mathcal{E}_{\gamma}(A),$ and thus of $\mathcal{E}(A)$, but it is also a
useful tool for understanding rank reversals for vectors in $\mathcal{E}%
_{\gamma}(A)$. Corollary \ref{Thintgea} also allows us to study the existence
of rank reversals for vectors in $\mathcal{E}(A)$. This will be elaborated in
Section \ref{srrev}.

\bigskip

We conclude this section with a result, that will be used later, providing an
interval for the product of the entries in $P_{A,\gamma}$ along a cycle.

\begin{lemma}
\label{ldecrea}Let $A\in\mathcal{PC}_{n}^{0}$ and $\gamma\in\Gamma(A)$. Let
$\beta(P_{A,\gamma})$ be the product of the entries in $P_{A,\gamma}$ along
the cycle $\beta=\mathcal{\beta}_{1}\cdots\mathcal{\beta}_{k}\mathcal{\beta
}_{1},$ with $2\leq k\leq n.$ Then,
\begin{equation}
\left(  \gamma(A)\right)  ^{k-1}\leq\beta(P_{A,\gamma})\leq\gamma(A).
\label{P1}%
\end{equation}

\end{lemma}

\begin{proof}
Let $w\in\mathcal{E}_{\gamma}(A)$. Then, by Theorem \ref{Thintgamma},
\begin{equation}
P_{A,\gamma}(\beta_{t},\beta_{t+1})\leq\frac{w_{\mathcal{\beta}_{t}}%
}{w_{\mathcal{\beta}_{t+1}}}\leq\frac{1}{P_{A,\gamma}(\mathcal{\beta}%
_{t+1},\mathcal{\beta}_{t})}=\frac{P_{A,\gamma}(\mathcal{\beta}_{t}%
,\mathcal{\beta}_{t+1})}{\gamma(A)}, \label{in12}%
\end{equation}
for any $t\in\{1,\ldots,k-1\}.$ From the first inequality in (\ref{in12}),
\[
\frac{w_{\mathcal{\beta}_{1}}}{w_{\mathcal{\beta}_{k}}}=\frac
{w_{\mathcal{\beta}_{1}}}{w_{\mathcal{\beta}_{2}}}\cdots\frac
{w_{\mathcal{\beta}_{k-1}}}{w_{\mathcal{\beta}_{k}}}\geq P_{A,\gamma
}(\mathcal{\beta}_{1},\mathcal{\beta}_{2})\cdots P_{A,\gamma}(\mathcal{\beta
}_{k-1},\mathcal{\beta}_{k}).
\]
On the other hand, by Theorem \ref{rattain}, there is a $w\in\mathcal{E}%
_{\gamma}(A)$ such that $\frac{w_{\mathcal{\beta}_{1}}}{w_{\mathcal{\beta}%
_{k}}}=P_{A,\gamma}(\mathcal{\beta}_{1},\mathcal{\beta}_{k})=\frac{\gamma
(A)}{P_{A,\gamma}(\mathcal{\beta}_{k},\mathcal{\beta}_{1})}.$ Therefore,%
\[
\frac{\gamma(A)}{P_{A,\gamma}(\mathcal{\beta}_{k},\mathcal{\beta}_{1})}\geq
P_{A,\gamma}(\mathcal{\beta}_{1},\mathcal{\beta}_{2})\cdots P_{A,\gamma
}(\mathcal{\beta}_{k-1},\mathcal{\beta}_{k}),
\]
implying the second inequality in (\ref{P1}).

From the second inequality in (\ref{in12}),
\[
\frac{w_{\mathcal{\beta}_{1}}}{w_{\mathcal{\beta}_{k}}}=\frac
{w_{\mathcal{\beta}_{1}}}{w_{\mathcal{\beta}_{2}}}\cdots\frac
{w_{\mathcal{\beta}_{k-1}}}{w_{\mathcal{\beta}_{k}}}\leq\frac{P_{A,\gamma
}(\mathcal{\beta}_{1},\mathcal{\beta}_{2})\cdots P_{A,\gamma}(\mathcal{\beta
}_{k-1},\mathcal{\beta}_{k})}{\left(  \gamma(A)\right)  ^{k-1}}.
\]
On the other hand, by Theorem \ref{rattain}, there is a $w\in\mathcal{E}%
_{\gamma}(A)$ such that $\frac{w_{\mathcal{\beta}_{1}}}{w_{\mathcal{\beta}%
_{k}}}=\frac{1}{P_{A,\gamma}(\mathcal{\beta}_{k},\mathcal{\beta}_{1})}.$
Therefore,%
\[
\frac{1}{P_{A,\gamma}(\mathcal{\beta}_{k},\mathcal{\beta}_{1})}\leq
\frac{P_{A,\gamma}(\mathcal{\beta}_{1},\mathcal{\beta}_{2})\cdots P_{A,\gamma
}(\mathcal{\beta}_{k-1},\mathcal{\beta}_{k})}{\left(  \gamma(A)\right)
^{k-1}},
\]
implying the first inequality in (\ref{P1}).
\end{proof}

\section{Existence of an efficient vector with a given order\label{Sorder}}

It may be important to know if there is an efficient vector $w$ for
$A\in\mathcal{PC}_{n}^{0}$ with a certain strict order. If $P$ is a
permutation matrix such that $Pw$ is decreasing, then, because of Lemma
\ref{lsim}, the question is equivalent to determining if there is a strictly
decreasing vector in $\mathcal{E}(PAP^{T})$ and, thus, in $\mathcal{E}%
_{\gamma}(PAP^{T})$ for some $\gamma\in\Gamma(PAP^{T}).$ This may be viewed as
a cone intersection problem, the two cones being the open cone of strictly
decreasing vectors in $\mathbb{R}_{+}^{n}$ and $\mathcal{E}_{\gamma}(A)$. Our
goal is to discuss the question in terms of the matrices $L_{A}$ or
$P_{A,\gamma}$.

\subsection{Sufficient conditions}

Here we give several conditions that guarantee the existence of a strictly
decreasing vector efficient for $A$.

\begin{theorem}
\label{ldecrLA}Let $A\in\mathcal{PC}_{n}^{0}$. If $L_{A}(j,j+1)>1$ for every
$j\in\{1,\ldots,n-1\}$, then every $w\in\mathcal{E}(A)$ is strictly
decreasing. If there is just one $j\in\{1,\ldots,n-1\}$ such that
$L_{A}(j,j+1)\leq1$ and, for such $j$, $L_{A}(j+1,j)<1$, then there is a
strictly decreasing vector $w\in\mathcal{E}(A)$.
\end{theorem}

\begin{proof}
If $L_{A}(j,j+1)>1$ for all $j\in\{1,\ldots,n-1\}$, by Corollary
\ref{Thintgea}, for any $w\in\mathcal{E}(A)$, $\frac{w_{j}}{w_{j+1}}>1$,
implying that $w$ is strictly decreasing. Suppose that there is exactly one
$j\in\{1,\ldots,n-1\}$ such that $L_{A}(j,j+1)\leq1,$ in which case, by
assumption, $L_{A}(j+1,j)<1$. By Theorem \ref{rattain}, there is a
$w\in\mathcal{E}(A)$ such that $\frac{w_{j+1}}{w_{j}}=L_{A}(j+1,j)<1$,
implying $\frac{w_{j}}{w_{j+1}}>1$. On the other hand, for $i\in
\{1,\ldots,n-1\}\backslash\{j\},$ $\frac{w_{i}}{w_{i+1}}\geq L_{A}(i,i+1)>1$,
in which the last inequality holds by assumption. Thus, $w$ is strictly decreasing.
\end{proof}

\bigskip

With a similar proof, now using Theorem \ref{Thintgamma} instead of Corollary
\ref{Thintgea}, we have the following sufficient conditions for the existence
of a strictly decreasing vector in $\mathcal{E}_{\gamma}(A)$.

\begin{theorem}
\label{ldecrgama}Let $A\in\mathcal{PC}_{n}^{0}$ and $\gamma\in\Gamma(A)$. If
$P_{A,\gamma}(j,j+1)>1$ for every $j\in\{1,\ldots,n-1\}$, then every
$w\in\mathcal{E}_{\gamma}(A)$ is strictly decreasing. If there is just one
$j\in\{1,\ldots,n-1\}$ such that $P_{A,\gamma}(j,j+1)\leq1$ and, for such $j$,
$P_{A,\gamma}(j+1,j)<1$, then there is a strictly decreasing vector
$w\in\mathcal{E}_{\gamma}(A)$.
\end{theorem}

\begin{theorem}
Let $A\in\mathcal{PC}_{n}^{0}$. Suppose that all the entries in $A$ along an
H-cycle $\gamma$ are $<1$. Then, $\gamma\in\Gamma(A)$ and there is a strictly
decreasing $w\in\mathcal{E}_{\gamma}(A)$.
\end{theorem}

\begin{proof}
Clearly, if all entries of $A$ along $\gamma$ are $<1$ then $\gamma\in
\Gamma(A)$ and all path products along $\gamma$ are $<1$ (that is, all
off-diagonal entries of $P_{A,\gamma}$ are $<1$). Thus, under the assumption,
we have $P_{A,\gamma}<J<P_{A,\gamma}^{(-T)}$ (except for the diagonal
positions). If $w$ is a strictly decreasing, sufficiently small, perturbation
of the $n$-vector of ones, then $P_{A,\gamma}\leq ww^{(-T)}\leq P_{A,\gamma
}^{(-T)}$. By Theorem \ref{Thintgamma}, $w\in\mathcal{E}_{\gamma}(A)$.
\end{proof}

\bigskip

The following result will be used in Section \ref{sbest}.

\begin{proposition}
\label{cdecrsuper}If all entries on the superdiagonal of $A\in\mathcal{PC}%
_{n}^{0}$ are $>1$, then there is a strictly decreasing vector $w\in
\mathcal{E}(A)$.
\end{proposition}

\begin{proof}
If $w\in\mathbb{R}_{+}^{n}$ is such that $ww^{(-T)}$ coincides with $A$ on the
superdiagonal, then $w$ is strictly decreasing and, by Remark
\ref{remspannintree}, $w\in\mathcal{E}(A)$.
\end{proof}

\bigskip

For $\gamma=12\cdots n1$, if $\gamma(A)<1$ then $\gamma\in\Gamma(A)$, and if
$\gamma(A)>1$ then $\operatorname*{rev}\gamma\in\Gamma(A).$ Consider the
vector $w$ mentioned in the above proof. If $\gamma(A)<1$, $w\in
\mathcal{E}_{\gamma}(A)$. If $\gamma(A)>1$, $w\in\mathcal{E}%
_{\operatorname*{rev}\gamma}(A)$. In either case, $w$ is an extreme of the set
$\mathcal{E}_{\gamma}(A)$ or $\mathcal{E}_{\operatorname*{rev}\gamma}(A)$. If
$A$ has superdiagonal entries $>1$ and $\gamma(A)<1$, the superdiagonal
entries of $P_{A,\gamma}$ are $>1$, so that the claim that $w$ is strictly
decreasing also follows from Theorem \ref{ldecrgama}. The following more
general fact is true.

\begin{corollary}
\label{ctr}Let $A\in\mathcal{PC}_{n}^{0}$ and suppose that $\gamma:=12\cdots
n1\in\Gamma(A)$. Then, all vectors in $\mathcal{E}_{\gamma}(A)$ are strictly
decreasing if and only if all entries on the superdiagonal of $A$ are $>1$.
\end{corollary}

\subsection{Necessary conditions}

Given $A\in\mathcal{PC}_{n}^{0}$ and $\gamma\in\Gamma(A)$, next we give
necessary conditions for the existence of a strictly decreasing vector in
$\mathcal{E}_{\gamma}(A)$. We denote by $\pi=\pi_{A,\gamma}$ the product of
the entries $a_{ij}$ along $\gamma$ such that either $a_{ij}\geq1$ or $i>j$
(that is, the product of the entries in $A$ along $\gamma$ with those above
the diagonal that are $<1$ replaced by $1$).

\begin{theorem}
\label{TT}Let $A\in\mathcal{PC}_{n}^{0}$ and $\gamma\in\Gamma(A)$. If there is
a strictly decreasing $w\in\mathcal{E}_{\gamma}(A)$, then $\pi_{A,\gamma}<1$
and all the entries of $P_{A,\gamma}$ below the diagonal are $<1$.
\end{theorem}

\begin{proof}
Suppose that there is a strictly decreasing vector $w\in\mathcal{E}_{\gamma
}(A)$. Since all entries of $ww^{(-T)}$ below the diagonal are $<1,$ by
Theorem \ref{Thintgamma}, all entries of $P_{A,\gamma}$ (and, thus, of $L_{A}%
$) below the diagonal are $<1$. (In particular, $a_{ij}<1$ for all edges
$(i,j)$ in $\gamma$, with $i>j,$ as $P_{A,\gamma}(i,j)=a_{ij}.$)

Let $\gamma=\gamma_{1}\cdots\gamma_{n}\gamma_{1}$ and suppose, without loss of
generality, that $\gamma_{1}=1$. Denote by $\pi^{\prime}$ the product of the
entries $a_{ij}$ along the path $\gamma_{1}\cdots\gamma_{n}$ such that either
$a_{ij}\geq1$ or $i>j.$ Since $\gamma_{n}>\gamma_{1}$, we have $\pi
=\pi^{\prime}a_{\gamma_{n}\gamma_{1}}.$

If all entries in $A$ along $\gamma$, and also above the diagonal, are $\geq
1$, then $\pi=\pi^{\prime}a_{\gamma_{n}\gamma_{1}}=\gamma(A)<1.$

Suppose that there is an entry in $A$ along $\gamma$, and also above the
diagonal, that is $<1$. By Theorem \ref{Thintgamma} (or the definition of
$\mathcal{E}_{\gamma}(A)$), we have%
\begin{equation}
\frac{w_{\gamma_{i}}}{w_{\gamma_{i+1}}}\geq a_{\gamma_{i}\gamma_{i+1}%
}\text{,\qquad}i\in\{1,\ldots,n-1\}, \label{gi}%
\end{equation}
and
\begin{equation}
\frac{w_{\gamma_{n}}}{w_{\gamma_{1}}}\geq a_{\gamma_{n}\gamma_{1}%
}\text{,\qquad that is,}\qquad\frac{w_{\gamma_{1}}}{w_{\gamma_{n}}}\leq
a_{\gamma_{1}\gamma_{n}}. \label{g1gn}%
\end{equation}
By (\ref{gi}) and using the fact that $w$ is decreasing, we have%
\[
\frac{w_{\gamma_{1}}}{w_{\gamma_{n}}}=\frac{w_{\gamma_{1}}}{w_{\gamma_{2}}%
}\cdots\frac{w_{\gamma_{n-1}}}{w_{\gamma_{n}}}>\pi^{\prime}\text{.}%
\]
Then, because of (\ref{g1gn}), $\pi^{\prime}<a_{\gamma_{1}\gamma_{n}},$ that
is, $\pi=\pi^{\prime}a_{\gamma_{n}\gamma_{1}}<1$.
\end{proof}

\bigskip

The first necessary condition in Theorem \ref{TT} does not imply the second.
In fact, we may have $\pi<1$ and some entry of $P_{A,\gamma}$ below the
diagonal $>1$, even for $n=3.$ In this case $\gamma=1321$, as if $\gamma=1231$
we have $P_{A,\gamma}(i,j)\leq\pi<1$ for $i>j$.

\begin{example}
Let $A\in\mathcal{PC}_{3}^{0}$ and $\gamma=1321$. Suppose that $a_{13}%
=\frac{1}{2},$ $a_{32}=\frac{1}{4}$ and $a_{21}=2$, so that $\Gamma
(A)=\{\gamma\}$. Then, $P_{A\,\gamma}(2,1)=a_{21}>1$ and $\pi=1.\frac{1}%
{4}.2<1.$
\end{example}

When $n>3$, we may have $\pi>1$ and all entries of $P_{A,\gamma}$ below the
diagonal $<1$. In this case, by Theorem \ref{TT} there is no strictly
decreasing vector in $\mathcal{E}_{\gamma}(A)$.

\begin{example}
Let $A\in\mathcal{PC}_{4}^{0}$ and $\gamma=12341$. Suppose that $a_{12}%
=a_{34}=\frac{1}{3},$ $a_{41}=\frac{1}{2}$ and $a_{23}=3.$ Then,
$\gamma(A)=\frac{1}{6}$, so that $\gamma\in\Gamma(A)$, and $\pi=1.3.1.\frac
{1}{2}=\frac{3}{2}>1.$ We have
\[
\mathcal{E}_{\gamma}(A)=\left\{  w\in\mathbb{R}_{+}^{4}:w_{1}\geq\frac{1}%
{3}w_{2}\geq w_{3}\geq\frac{1}{3}w_{4}\geq\frac{1}{6}w_{1}\right\}  .
\]
Also,
\[
P_{A,\gamma}=\left[
\begin{array}
[c]{cccc}%
1 & \frac{1}{3} & 1 & \frac{1}{3}\\
\frac{1}{2} & 1 & 3 & 1\\
\frac{1}{6} & \frac{1}{18} & 1 & \frac{1}{3}\\
\frac{1}{2} & \frac{1}{6} & \frac{1}{2} & 1
\end{array}
\right]  ,
\]
which has all entries below the diagonal $<1$. If $w$ were a strictly
decreasing vector in $\mathcal{E}_{\gamma}(A),$ using Theorem \ref{Thintgamma}
and the fact that $\frac{w_{3}}{w_{4}}>1$, we would have
\[
\frac{w_{2}}{w_{1}}=\frac{w_{2}}{w_{3}}\frac{w_{3}}{w_{4}}\frac{w_{4}}{w_{1}%
}>3.1.\frac{1}{2}=\frac{3}{2}>1,
\]
a contradiction. Thus, there is no strictly decreasing vector in
$\mathcal{E}_{\gamma}(A)$, as expected.
\end{example}

When $n>4$, the two necessary conditions in Theorem \ref{TT} are not sufficient.

\begin{example}
Let $A\in\mathcal{PC}_{5}^{0}$ and $\gamma=152341$. Suppose that $a_{15}=2,$
$a_{52}=\frac{1}{2}$, $a_{23}=a_{41}=\frac{1}{4}$ and $a_{34}=3.$ Then,
$\gamma(A)=\frac{3}{16}$, so that $\gamma\in\Gamma(A)$, and $\pi=2.\frac{1}%
{2}.1.3.\frac{1}{4}=\frac{3}{4}<1.$ We have
\[
\mathcal{E}_{\gamma}(A)=\left\{  w\in\mathbb{R}_{+}^{5}:w_{1}\geq2w_{5}\geq
w_{2}\geq\frac{1}{4}w_{3}\geq\frac{3}{4}w_{4}\geq\frac{3}{16}w_{1}\right\}  .
\]
Also,
\[
P_{A,\gamma}=\left[
\begin{array}
[c]{ccccc}%
1 & 1 & \frac{1}{4} & \frac{3}{4} & 2\\
\frac{3}{16} & 1 & \frac{1}{4} & \frac{3}{4} & \frac{3}{8}\\
\frac{3}{4} & \frac{3}{4} & 1 & 3 & \frac{3}{2}\\
\frac{1}{4} & \frac{1}{4} & \frac{1}{16} & 1 & \frac{1}{2}\\
\frac{3}{32} & \frac{1}{2} & \frac{1}{8} & \frac{3}{8} & 1
\end{array}
\right]  ,
\]
which has all entries below the diagonal $<1$. If $w$ were a strictly
decreasing vector in $\mathcal{E}_{\gamma}(A),$ using Theorem \ref{Thintgamma}
and the fact that $\frac{w_{2}}{w_{3}}>1$, we would have
\[
\frac{w_{5}}{w_{4}}=\frac{w_{5}}{w_{2}}\frac{w_{2}}{w_{3}}\frac{w_{3}}{w_{4}%
}\geq\frac{1}{2}.1.3=\frac{3}{2}>1,
\]
a contradiction. Thus, there is no strictly decreasing vector in
$\mathcal{E}_{\gamma}(A)$.
\end{example}

\subsection{The case $n=3$}

When $n=3$, the condition that all entries of $P_{A,\gamma}$ below the
diagonal are $<1$ implies that $\pi<1$. Moreover, it is sufficient for the
existence of a decreasing vector in $\mathcal{E}_{\gamma}(A)$. Recall that,
for $A\in\mathcal{PC}_{3}^{0},$ we have $\#\Gamma(A)=1$, so that if
$\Gamma(A)=\{\gamma\}$ then $\mathcal{E}(A)=\mathcal{E}_{\gamma}(A)$ and
$P_{A,\gamma}=L_{A}$.

\begin{theorem}
\label{n3}Let $A\in\mathcal{PC}_{3}^{0}$ and $\Gamma(A)=\{\gamma\}$. There is
a strictly decreasing $w\in\mathcal{E}_{\gamma}(A)$ if and only if all entries
of $P_{A,\gamma}$ below the diagonal are $<1$.
\end{theorem}

\begin{proof}
From Theorem \ref{TT}, if there is a strictly decreasing $w\in\mathcal{E}%
_{\gamma}(A)$ then all entries of $P_{A,\gamma}$ below the diagonal are $<1.$

To see the converse, suppose that all entries of $P_{A,\gamma}$ below the
diagonal are $<1$. If there is at most one $i\in\{1,2\}$ such that
$P_{A,\gamma}(i,i+1)\leq1$, the claim follows from Theorem \ref{ldecrgama}.

Suppose now that $P_{A,\gamma}(1,2)\leq1$ and $P_{A,\gamma}(2,3)\leq1$. By
Lemma \ref{ldecrea}, we have $P_{A,\gamma}(2,1)P_{A,\gamma}(3,2)\leq
P_{A,\gamma}(3,1)(=\frac{\gamma(A)}{P_{A,\gamma}(1,3)}).$ Then, since all
entries of $P_{A,\gamma}$ below the diagonal are $<1$, write $P_{A,\gamma
}(3,1)=x_{1}x_{2}$, with $P_{A,\gamma}(2,1)\leq x_{1}<1$ and $P_{A,\gamma
}(3,2)\leq x_{2}<1.$ Let $w\in\mathbb{R}_{+}^{3}$ be such that $\frac{w_{2}%
}{w_{1}}=x_{1}$ and $\frac{w_{3}}{w_{2}}=x_{2}$. Clearly, $w$ is strictly
decreasing. Also, noticing that $\frac{1}{P_{A,\gamma}(1,2)}\geq1$, $\frac
{1}{P_{A,\gamma}(2,3)}\geq1$ and $P_{A,\gamma}(3,1)=\frac{w_{3}}{w_{1}}%
<\frac{1}{P_{A,\gamma}(1,3)},$ by Theorem \ref{Thintgamma}, $w\in
\mathcal{E}_{\gamma}(A)$.
\end{proof}

\begin{example}
Let%
\[
A=\left[
\begin{array}
[c]{ccc}%
1 & \frac{1}{2} & 4\\
2 & 1 & \frac{1}{2}\\
\frac{1}{4} & 2 & 1
\end{array}
\right]  \in\mathcal{PC}_{3}^{0}.
\]
Then $\Gamma(A)=\{\gamma\}$ with $\gamma=1231,$ as $\gamma(A)=\frac{1}{16}.$
We have%
\[
P_{A,\gamma}=\left[
\begin{array}
[c]{ccc}%
1 & \frac{1}{2} & \frac{1}{4}\\
\frac{1}{8} & 1 & \frac{1}{2}\\
\frac{1}{4} & \frac{1}{8} & 1
\end{array}
\right]  \in\mathcal{PC}_{3}^{0}.
\]
Since $P_{A,\gamma}(3,1)=\frac{1}{4}$, we may write $P_{A,\gamma}%
(3,1)=x_{1}x_{2}=\frac{1}{2}.\frac{1}{2}$. Let $w\in\mathbb{R}_{+}^{3}$ be
such that $\frac{w_{2}}{w_{1}}=\frac{1}{2}$ and $\frac{w_{3}}{w_{2}}=\frac
{1}{2}$; projectively, we have $w=\left[
\begin{array}
[c]{ccc}%
4 & 2 & 1
\end{array}
\right]  ^{T}$. Then, $w$ is strictly decreasing and, by Theorem
\ref{Thintgamma} (or Theorem \ref{blanq}), $w\in\mathcal{E}_{\gamma
}(A)=\mathcal{E}(A)$.
\end{example}

\section{Order ambiguity in $\mathcal{E}(A)$\label{s3}}

An $i,j$ order ambiguity in $\mathcal{E}(A)$ (resp. $\mathcal{E}_{\gamma}(A)$)
means that there are $w,v\in\mathcal{E}(A)$ (resp. $\mathcal{E}_{\gamma}(A)$)
such that $\frac{w_{i}}{w_{j}}<1$ and $\frac{v_{i}}{v_{j}}>1$. In this case,
not all vectors in $\mathcal{E}(A)$ (resp. $\mathcal{E}_{\gamma}(A)$) give the
same ordering of alternatives $i$ and $j$. The results in Section
\ref{sbackground} allow us to characterize order ambiguity in $\mathcal{E}(A)$
and $\mathcal{E}_{\gamma}(A)$ in terms of $L_{A}$ and $P_{A,\gamma}$, respectively.

\begin{theorem}
\label{Porderambgamma}Let $A\in\mathcal{PC}_{n}^{0}$, $\gamma\in\Gamma(A)$,
and $i,j\in\{1,\ldots,n\}$, $i\neq j.$ There is an $i,j$ order ambiguity in
$\mathcal{E}_{\gamma}(A)$ (resp. $\mathcal{E}(A)$) if and only if both
$P_{A,\gamma}(i,j)<1$ and $P_{A,\gamma}(j,i)<1$ (resp. both $L_{A}(i,j)<1$ and
$L_{A}(j,i)<1$).
\end{theorem}

\begin{proof}
By Theorems \ref{Thintgamma} and \ref{rattain}, there is an $i,j$ order
ambiguity in $\mathcal{E}_{\gamma}(A)$ if and only if $P_{A,\gamma}(i,j)<1$
and $\frac{1}{P_{A,\gamma}(j,i)}>1$. Using Corollary \ref{Thintgea} instead of
Theorem \ref{Thintgamma}, we have the corresponding result for $\mathcal{E}%
(A)$.
\end{proof}

\begin{corollary}
Let $A\in\mathcal{PC}_{n}^{0}$, $\gamma\in\Gamma(A)$, and $i,j\in
\{1,\ldots,n\}$, $i\neq j.$ In $\mathcal{E}(A)$ or $\mathcal{E}_{\gamma}(A)$,
there is an $i,j$ order ambiguity if and only if there is a $j,i$ order ambiguity.
\end{corollary}

\bigskip

Notice that, if $(i,j)$ is an edge in $\gamma\in\Gamma(A)$ and $a_{ij}%
(=P_{A,\gamma}(i,j))\geq1$, there is no $i,j$ order ambiguity in
$\mathcal{E}_{\gamma}(A)$.

\bigskip

It may happen that there is no $i,j$ order ambiguity in $\mathcal{E}_{\gamma
}(A)$ for any $i,j$, $i\neq j$ (see Section \ref{sUORR}).

\begin{example}
\label{Exnoorderamb}Continue Examples \ref{exbound0} and \ref{exbound02}.
There is no order ambiguity in $\mathcal{E}_{\alpha}(A)$. Also, there is no
order ambiguity in either $\mathcal{E}_{\beta}(A)$ or $\mathcal{E}(A)$, except
at positions $2,3$ and $3,2$, where there is order ambiguity.
\end{example}

Of course, an $i,j$ order ambiguity in $\mathcal{E}_{\gamma}(A),$ for some
$\gamma\in\Gamma(A)$, implies an $i,j$ order ambiguity in $\mathcal{E}(A)$. As
one might guess, the converse is not true. However, if for all $\gamma
\in\Gamma(A)$ there is an $i,j$ order unambiguity in $\mathcal{E}_{\gamma}(A)$
of the same type, that is, $P_{A,\gamma}(i,j)\geq1$ for all $\gamma\in
\Gamma(A),$ or $P_{A,\gamma}(j,i)\geq1$ for all $\gamma\in\Gamma(A),$ then
there is an $i,j$ order unambiguity in $\mathcal{E}(A).$

\bigskip

For $A\in\mathcal{PC}_{n}^{0}$ and $\gamma\in\Gamma(A)$, we next characterize
when there is an $i,j$ order ambiguity in $\mathcal{E}_{\gamma}(A)$ (and,
thus, in $\mathcal{E}(A)$) for all $i\neq j$.

\begin{theorem}
\label{thallorderamb}Let $A\in\mathcal{PC}_{n}^{0}$ and $\gamma\in\Gamma(A)$.
There is an $i,j$ order ambiguity in $\mathcal{E}_{\gamma}(A)$ for all
$i,j\in\{1,\ldots,n\}$, $i\neq j$, if and only if all entries in $A$ along
$\gamma$ are $<1$.
\end{theorem}

\begin{proof}
By Theorem \ref{Porderambgamma}, there is an $i,j$ order ambiguity in
$\mathcal{E}_{\gamma}(A)$ if and only if $P_{A,\gamma}(i,j)<1$ and
$P_{A,\gamma}(j,i)<1$. Clearly, if all entries in $A$ along $\gamma$ are $<1$
then $P_{A,\gamma}(i,j)<1$ for all $i,j\in\{1,\ldots,n\}$, $i\neq j$. The
converse follows because, if $(i,j)$ is an edge in $\gamma$, then
$P_{A,\gamma}(i,j)=a_{ij}$.
\end{proof}

\begin{example}
\label{Exnoamb}Let%
\[
A=\left[
\begin{array}
[c]{cccc}%
1 & \frac{1}{2} & 2 & 8\\
2 & 1 & \frac{1}{2} & 2\\
\frac{1}{2} & 2 & 1 & \frac{1}{2}\\
\frac{1}{8} & \frac{1}{2} & 2 & 1
\end{array}
\right]  \in\mathcal{PC}_{4}^{(1)}.
\]
We have $\Gamma(A)=\{\alpha\},$ with $\alpha=12341$, and $\alpha(A)=\frac
{1}{64}$. Thus,
\[
\mathcal{E}(A)=\mathcal{E}_{\alpha}(A)=\left\{  w\in\mathbb{R}_{+}^{4}%
:w_{1}\geq\frac{1}{2}w_{2}\geq\frac{1}{4}w_{3}\geq\frac{1}{8}w_{4}\geq\frac
{1}{64}w_{1}\right\}  .
\]
For $W=ww^{(-T)},$ with $w\in\mathcal{E}(A)$, we have
\[
L_{A}=P_{A,\alpha}=\left[
\begin{array}
[c]{cccc}%
1 & \frac{1}{2} & \frac{1}{4} & \frac{1}{8}\\
\frac{1}{32} & 1 & \frac{1}{2} & \frac{1}{4}\\
\frac{1}{16} & \frac{1}{32} & 1 & \frac{1}{2}\\
\frac{1}{8} & \frac{1}{16} & \frac{1}{32} & 1
\end{array}
\right]  \leq W\leq\left[
\begin{array}
[c]{cccc}%
1 & 32 & 16 & 8\\
2 & 1 & 32 & 16\\
4 & 2 & 1 & 32\\
8 & 4 & 2 & 1
\end{array}
\right]  =P_{A,\alpha}^{(-T)}=L_{A}^{(-T)}.
\]
There is an $i,j$ order ambiguity in $\mathcal{E}(A)=\mathcal{E}_{\alpha}(A)$,
for all $i\neq j.$
\end{example}

\section{Rank reversals in $\mathcal{E}(A)$\label{srrev}}

In this section we focus on how rank reversals occur. Given $A\in
\mathcal{PC}_{n}^{0}$ and $i,j\in\{1,\ldots,n\},$ $i\neq j$, there is a
$w\in\mathcal{E}(A)$ such that $\frac{w_{i}}{w_{j}}=a_{ij}$. Such a $w$ has no
$i,j$ rank reversal. One example is the $j$th column of $A$, which is in
$\mathcal{E}(A)$ \cite{FJ1} and has no rank reversals for $A$ at positions in
column and row $j$. However, rank reversals at other positions may be
unavoidable and we consider here how they occur.

\bigskip

Clearly, if $w\in\mathbb{R}_{+}^{n}$ has an $i,j$ rank reversal, it also has a
$j,i$ rank reversal. The tally of rank reversals for a vector is invariant
under permutation similarity.

\begin{proposition}
\label{Prankpermut}Let $A\in\mathcal{PC}_{n}^{0}$, $P$ be an $n$-by-$n$
permutation matrix and $w\in\mathbb{R}_{+}^{n}$. Then, the number of rank
reversals for $A$ at $w$ is the same as for $PAP^{T}$ at $Pw$.
\end{proposition}

For simplicity, in this section we focus on matrices $A\in\mathcal{PC}%
_{n}^{(1)}$.

\subsection{The relationship with order ambiguity}

Now we record several observations about rank reversals, and their
relationship with order ambiguity, that are valid because of Theorems
\ref{Thintgamma} and \ref{rattain} and Corollary \ref{Thintgea}.

\bigskip

\begin{theorem}
\label{propanbrev}Let $A\in\mathcal{PC}_{n}^{(1)}$ and $i,j\in\{1,\ldots,n\},$
$i\neq j$. There is an $i,j$ order ambiguity in $\mathcal{E}(A)$ if and only
if there is a $w\in\mathcal{E}(A)$ with an $i,j$ rank reversal.
\end{theorem}

\begin{proof}
If there is an $i,j$ order ambiguity in $\mathcal{E}(A)$, then, by Theorem
\ref{Porderambgamma}, $L_{A}(i,j)<1$ and $\frac{1}{L_{A}(j,i)}>1.$ If
$a_{ij}>1$, the vector $w\in\mathcal{E}(A)$ such that $\frac{w_{i}}{w_{j}%
}=L_{A}(i,j)$ (whose existence is guaranteed by Theorem \ref{rattain}) has an
$i,j$ rank reversal. If $a_{ij}<1$, the vector $w\in\mathcal{E}(A)$ such that
$\frac{w_{i}}{w_{j}}=\frac{1}{L_{A}(j,i)}$ has an $i,j$ rank reversal. The
converse follows from the two claims about matrix intervals in Corollary
\ref{Thintgea}.
\end{proof}

\begin{example}
\label{exbound03}Continue Examples \ref{exbound0}, \ref{exbound02} and
\ref{Exnoorderamb}. There is a vector in $\mathcal{E}(A)$ with a rank reversal
at positions $2,3$ and $3,2$. For each of the other positions, there is no
vector in $\mathcal{E}(A)$ with a rank reversal.
\end{example}

When $\mathcal{E}_{\gamma}(A)$ is considered, it may happen that $a_{ij}$ is
not in the interval $\left[  P_{A,\gamma}(i,j),\frac{1}{P_{A,\gamma}%
(j,i)}\right]  $. In this situation, $(i,j)$ and $(j,i)$ are not edges in
$\gamma$, as otherwise $P_{A,\gamma}(i,j)=a_{ij}$ or $\frac{1}{P_{A,\gamma
}(j,i)}=a_{ij}$.

\begin{theorem}
\label{Porderreversalgamma}Let $A\in\mathcal{PC}_{n}^{(1)}$, $\gamma\in
\Gamma(A)$ and $i,j\in\{1,\ldots,n\},$ $i\neq j$. If there is an $i,j$ order
ambiguity in $\mathcal{E}_{\gamma}(A),$ then there is a $w\in\mathcal{E}%
_{\gamma}(A)$ with an $i,j$ rank reversal. Conversely, if there is a
$w\in\mathcal{E}_{\gamma}(A)$ with an $i,j$ rank reversal then either there is
an $i,j$ order ambiguity in $\mathcal{E}_{\gamma}(A),$ or there is an $i,j$
rank reversal for all vectors $w\in\mathcal{E}_{\gamma}(A)$ with $w_{i}\neq
w_{j}.$ In the latter case, $a_{ij}\notin\left[  P_{A,\gamma}(i,j),\frac
{1}{P_{A,\gamma}(j,i)}\right]  $ and, thus, $(i,j)$ and $(j,i)$ are not edges
in $\gamma$.
\end{theorem}

\begin{proof}
The first claim follows by the same argument as in the proof of the
corresponding claim in Theorem \ref{propanbrev}, replacing $\mathcal{E}(A)$ by
$\mathcal{E}_{\gamma}(A)$ and $L_{A}$ by $P_{A,\gamma}$.

Regarding the second claim, suppose that there is a $w\in\mathcal{E}_{\gamma
}(A)$ with an $i,j$ rank reversal. Then, because of Theorem \ref{Thintgamma},
either $P_{A,\gamma}(i,j)<1$ and $a_{ij}>1,$ or $\frac{1}{P_{A,\gamma}%
(j,i)}>1$ and $a_{ij}<1$. If $a_{ij}\in\left[  P_{A,\gamma}(i,j),\frac
{1}{P_{A,\gamma}(j,i)}\right]  $, then there is an $i,j$ order ambiguity in
$\mathcal{E}_{\gamma}(A)$. If $a_{ij}\notin\left[  P_{A,\gamma}(i,j),\frac
{1}{P_{A,\gamma}(j,i)}\right]  $ and there is no $i,j$ order ambiguity in
$\mathcal{E}_{\gamma}(A)$, then either $a_{ij}<1$ and $P_{A,\gamma}%
(i,j)\geq1,$ or $a_{ij}>1$ and $\frac{1}{P_{A,\gamma}(j,i)}\leq1.$ By Theorem
\ref{Thintgamma}, there is an $i,j$ rank reversal for all vectors
$w\in\mathcal{E}_{\gamma}(A)$, with $w_{i}\neq w_{j}$ (if there is a $w$ with
$w_{i}=w_{j}$ then $P_{A,\gamma}(i,j)=1$ or $P_{A,\gamma}(j,i)=1$)$.$
\end{proof}

\begin{example}
\label{exrr}Continue Examples \ref{exbound0}, \ref{exbound02},
\ref{Exnoorderamb} and \ref{exbound03}. Since in $\mathcal{E}_{\beta}(A)$
there is order ambiguity just at positions $2,3$ and $3,2$, and $P_{A,\beta
}\leq A\leq P_{A,\beta}^{(-T)}$, the only possible rank reversals for vectors
in $\mathcal{E}_{\beta}(A)$ occur at positions $2,3$ and $3,2$. (There is a
vector $w\in\mathcal{E}_{\beta}(A)$ with rank reversals at these positions and
also a vector $v\in\mathcal{E}_{\beta}(A)$ with no rank reversals at the same
positions.) In $\mathcal{E}_{\alpha}(A)$ there are no vectors with rank
reversals. This is because there is no order ambiguity in $\mathcal{E}%
_{\alpha}(A)$, and $P_{A,\alpha}\leq A\leq P_{A,\alpha}^{(-T)}$, except in
positions $1,3$ $(3,1)$ and $2,4$ ($4,2$). However, in such positions $i,j$ we
still have that $P_{A,\gamma}(i,j)>1$ if $a_{ij}>1$ and $\frac{1}{P_{A,\gamma
}(j,i)}<1$ if $a_{ij}<1$.

Consider the matrix $B\in\mathcal{PC}_{4}^{(1)}$ obtained from the matrix $A$
in Example \ref{exbound0} by replacing the entries in positions $1,3$ and
$3,1$ by $\frac{1}{2}$ and $2,$ respectively. Then $\alpha\in\Gamma(B).$
Moreover, $\mathcal{E}_{\alpha}(B)=\mathcal{E}_{\alpha}(A)$ and $P_{B,\alpha
}=P_{A,\alpha}$, since $A$ and $B$ have the same entries along $\alpha$.
However, all vectors in $\mathcal{E}_{\alpha}(B)$ have a rank reversal for $B$
at positions $1,3$ and $3,1$.
\end{example}

\subsection{Rank reversals in $\mathcal{E}_{\gamma}(A)$ along $\gamma$}

\color{black}Since $\mathcal{E}_{\gamma}(A)$ just depends on the entries in
$A$ along $\gamma$, if neither $(i,j)$ nor $(j,i)$ is an edge in $\gamma$, the
existence of an $i,j$ rank reversal at a given $w\in\mathcal{E}_{\gamma}(A)$
cannot be predicted from $\mathcal{E}_{\gamma}(A)$, as indicated by Theorem
\ref{Porderreversalgamma} and illustrated in Example \ref{exrr}. However, in
positions along $\gamma$ (and, thus, $\operatorname*{rev}\gamma$), some
conclusions can be drawn.

\begin{proposition}
\label{cga}Let $A\in\mathcal{PC}_{n}^{(1)}$, $\gamma\in\Gamma(A)$ and let
$(i,j)$ be an edge in $\gamma$. If $a_{ij}>1$ or $\gamma(A)a_{ji}\geq1$, there
is no $w\in\mathcal{E}_{\gamma}(A)$ with an $i,j$ rank reversal.
\end{proposition}

\begin{proof}
Since $(i,j)$ is an edge in $\gamma$, $P_{A,\gamma}(i,j)=a_{ij},$ implying
that $a_{ij}\in\lbrack P_{A,\gamma}(i,j),\frac{1}{P_{A,\gamma}(j,i)}]$. Since
$P_{A,\gamma}(j,i)=\gamma(A)a_{ji}$, under the assumptions, Theorem
\ref{Porderambgamma} ensures that there is no $i,j$ order ambiguity in
$\mathcal{E}_{\gamma}(A)$. Now the claim follows from Theorem
\ref{Porderreversalgamma}.
\end{proof}

\bigskip

Though it may happen that a vector in $\mathcal{E}_{\gamma}(A)$, with
$\gamma\in\Gamma(A)$, has a rank reversal in some positions along $\gamma$, it
cannot have a rank reversal in all positions along $\gamma$.

\begin{theorem}
Let $A\in\mathcal{PC}_{n}^{(1)}$, $\gamma\in\Gamma(A)$ and $w\in
\mathcal{E}_{\gamma}(A)$. Then, there is an edge $(i,j)$ in $\gamma$ such that
$w$ has no rank reversal at $i,j$.
\end{theorem}

\begin{proof}
If there is an entry $a_{ij}>1$ with $(i,j)$ an edge in $\gamma$, then, by
Proposition \ref{cga}, $w$ has no $i,j$ rank reversal. Suppose that all
entries $a_{ij}$ along $\gamma$ are $<1.\ $Since
\begin{equation}
\frac{w_{\gamma_{n}}}{w_{\gamma_{1}}}=\frac{w_{\gamma_{n}}}{w_{\gamma_{n-1}}%
}\cdots\frac{w_{\gamma_{2}}}{w_{\gamma_{1}}}\text{,} \label{qq}%
\end{equation}
if $\frac{w_{\gamma_{t}}}{w_{\gamma_{t+1}}}>1$ for all $1\leq t\leq n-1,$ then
$\frac{w_{\gamma_{n}}}{w_{\gamma_{1}}}<1$. Thus, $w$ has no rank reversal at
position $\gamma_{n},\gamma_{1}$.
\end{proof}

\bigskip

For any path with $n$ vertices in an H-cycle $\gamma\in\Gamma(A),$ there is a
(projectively unique) vector $w\in\mathcal{E}_{\gamma}(A)$ such that
$ww^{(-T)}$ coincides with $A$ along the path. So, along the H-cycle, this
vector gives at most one rank reversal. In some circumstances, a vector in
$\mathcal{E}_{\gamma}(A)$ may have no rank reversals along $\gamma$.

\begin{theorem}
\label{t23}Let $A\in\mathcal{PC}_{n}^{(1)}$ and $\gamma\in\Gamma(A)$. If there
is an entry in $A$ along $\gamma$ that is $>1$, then there is a strict
$w\in\mathcal{E}_{\gamma}(A)$ with no rank reversals along $\gamma$.
Conversely, if there is a nonconstant $w\in\mathcal{E}_{\gamma}(A)$ with no
rank reversals along $\gamma$, then there is an entry in $A$ along $\gamma$
that is $>1$.
\end{theorem}

\begin{proof}
Suppose that there is an entry in $A$ along $\gamma$ that is $>1$, in which
case we assume, without loss of generality, that $a_{\gamma_{n}\gamma_{1}}>1$.
Let $w$ be such that $\frac{w_{\gamma_{i}}}{w_{\gamma_{i+1}}}=a_{\gamma
_{i}\gamma_{i+1}},$ for $1\leq i\leq n-1.$ Then $w\in\mathcal{E}_{\gamma}(A)$.
Because of (\ref{qq}), we have%
\[
\frac{w_{\gamma_{n}}}{w_{\gamma_{1}}}=\frac{1}{a_{\gamma_{1}\gamma_{2}}\cdots
a_{\gamma_{n-1}\gamma_{n}}}=\frac{a_{\gamma_{n}\gamma_{1}}}{\gamma
(A)}>a_{\gamma_{n}\gamma_{1}}>1,
\]
so that $w$ has no rank reversals along $\gamma.$ If $w$ happens to be strict,
we are done. If $w$ is not strict, since $\mathcal{E}_{\gamma}(A)$ is convex
with positive density, there is a sufficiently small perturbation of $w$ that
is strict, is in $\mathcal{E}_{\gamma}(A)$ and still has no rank reversals
along $\gamma$.

For the converse, suppose that $w\in\mathcal{E}_{\gamma}(A)$ is nonconstant
and has no rank reversals along $\gamma$. Suppose that all entries in $A$
along the path $\gamma_{1}\cdots\gamma_{n}$ are $<1$. Then, $\frac
{w_{\gamma_{t}}}{w_{\gamma_{t+1}}}\leq1$ for all $1\leq t\leq n-1$, with at
least one inequality strict. Hence, because of (\ref{qq}), we have
$\frac{w_{\gamma_{n}}}{w_{\gamma_{1}}}>1$. Since $A$ has no rank reversals,
this implies that $a_{\gamma_{n}\gamma_{1}}>1$.
\end{proof}

\bigskip

Theorem \ref{t23} might be compared with Theorem 3.1 in \cite{FJ3}.

\section{Weight vectors with the fewest rank reversals\label{srankmin}}

The goal of this section is to show that, for any $A\in\mathcal{PC}_{n}^{(1)}%
$, there is a vector efficient for $A$ that gives the fewest rank reversals
among all strict vectors in $\mathbb{R}_{+}^{n},$ and to describe how to
obtain such a vector. Next we introduce some ideas to this end.

\subsection{s-minimizers}

We denote by $s(A)$ the minimum number of entries $<1$ above the diagonal,
among all matrices in the permutation similarity class of $A\in\mathcal{PC}%
_{n}^{(1)}$. A matrix $B$ permutationally similar to $A$ with $s(A)$ entries
$<1$ above the diagonal is called an $s$\emph{-minimizer} (in the permutation
similarity class of $A$). Clearly, $s(A)\leq\frac{n^{2}-n}{4}$, and $s(A)=0$
if and only if all entries above the diagonal of an $s$\emph{-minimizer} of
$A$ are $>1$. If $B$ is an $s$-minimizer, then the submatrix of $B$ obtained
by deleting the first row and column, or the last row and column, is also an
$s$-minimizer in the corresponding class. Notice also that, if $A$ has a row
with all off-diagonal entries $>1$, then the first row of any s-minimizer of
$A$ will be that row.

\bigskip

The problem of determining $s(A)$ and an $s$-minimizer for $A$ has an
equivalent formulation for tournaments. It is determination of the fewest
edges whose reversal gives a transitive tournament. This number of edges for
the tournament $T(A)$ is the same as our $s(A)$. This problem for tournaments
is known to be NP-hard \cite{Alon}. Algorithms and heuristic methods to solve
it are discussed, for example, in \cite{Ali}. For more recent literature see
\cite{Chen} and the references therein. Here, we are not concerned with these
difficulties due to size. In the reciprocal matrix literature the problems
tend not to be large.

\bigskip

An important fact is the following. The analog for tournaments has been
noticed in \cite{Ali}.

\begin{lemma}
\label{PP1}A matrix $B\in\mathcal{PC}_{n}^{(1)}$ that is an s-minimizer has
all entries on the superdiagonal $>1.$
\end{lemma}

\begin{proof}
Suppose that the entry in position $i,i+1$ of $B$ is $<1.$ Then, interchanging
rows and columns $i$ and $i+1$ decreases the number of entries $<1$ above the
diagonal, which contradicts the assumption that $B$ is an s-minimizer. Note
that, after permutation, all entries above the diagonal remain the same except
entry $b_{i,i+1}$ that is replaced by $b_{i+1,i}.$
\end{proof}

\bigskip A consequence of Lemma \ref{PP1} is that any matrix in $\mathcal{PC}%
_{n}^{(1)}$ is permutationally similar to a matrix with superdiagonal entries
$>1$. In the context of tournaments, this means that any tournament has a
Hamiltonian path. This is precisely R\'{e}dei's Theorem for tournaments
\cite{Moon}. If we consider the tournament $T(A)$, with $A\in\mathcal{PC}%
_{n}^{(1)}$, the same permutation that gives an s-minimizer for $A$, applied
to $T(A)$, gives a tournament matrix with superdiagonal entries $>1$.

\bigskip

If $A\in\mathcal{PC}_{n}^{(1)}$, with $3\leq n\leq4,$ then $s(A)\leq1$. Using
Lemma \ref{PP1}, the case $n=3$ is clear. Suppose that $n=4$ and suppose,
without loss of generality, that $A$ has exactly $s(A)$ entries $<1$ above the
diagonal. By Lemma \ref{PP1}, all entries on the superdiagonal of $A$ are
$>1$, implying that $s(A)\leq3$. Suppose that $2\leq s(A)\leq3.$ Then, $T(A)$
is one of the following matrices:
\[
\left[
\begin{array}
[c]{cccc}%
0 & 1 & 0 & 0\\
0 & 0 & 1 & 0\\
1 & 0 & 0 & 1\\
1 & 1 & 0 & 0
\end{array}
\right]  ,\text{ }\left[
\begin{array}
[c]{cccc}%
0 & 1 & 0 & 1\\
0 & 0 & 1 & 0\\
1 & 0 & 0 & 1\\
0 & 1 & 0 & 0
\end{array}
\right]  ,\text{ }\left[
\begin{array}
[c]{cccc}%
0 & 1 & 1 & 0\\
0 & 0 & 1 & 0\\
0 & 0 & 0 & 1\\
1 & 1 & 0 & 0
\end{array}
\right]  ,\text{ }\left[
\begin{array}
[c]{cccc}%
0 & 1 & 0 & 0\\
0 & 0 & 1 & 1\\
1 & 0 & 0 & 1\\
1 & 0 & 0 & 0
\end{array}
\right]  .
\]
For each case, there is a permutation matrix $P$ such that $PAP^{T}$ has
exactly one entry $<1$ above the diagonal, a contradiction. These permutation
matrices $P$ are%
\[
\left[
\begin{array}
[c]{cccc}%
0 & 0 & 1 & 0\\
0 & 0 & 0 & 1\\
1 & 0 & 0 & 0\\
0 & 1 & 0 & 0
\end{array}
\right]  ,\text{ }\left[
\begin{array}
[c]{cccc}%
0 & 0 & 1 & 0\\
1 & 0 & 0 & 0\\
0 & 0 & 0 & 1\\
0 & 1 & 0 & 0
\end{array}
\right]  ,\text{ }\left[
\begin{array}
[c]{cccc}%
0 & 0 & 0 & 1\\
1 & 0 & 0 & 0\\
0 & 1 & 0 & 0\\
0 & 0 & 1 & 0
\end{array}
\right]  ,\text{ }\left[
\begin{array}
[c]{cccc}%
0 & 1 & 0 & 0\\
0 & 0 & 1 & 0\\
0 & 0 & 0 & 1\\
1 & 0 & 0 & 0
\end{array}
\right]  ,
\]
for each case.

\bigskip

Computationally, we have verified the following. If $A\in\mathcal{PC}%
_{5}^{(1)}$, then $s(A)$ can be $0,1,2$ or $3$. For example, any matrix
$A\in\mathcal{PC}_{5}^{(1)}$ such that
\[
T(A)=\left[
\begin{array}
[c]{ccccc}%
0 & 1 & 1 & 0 & 0\\
0 & 0 & 1 & 1 & 0\\
0 & 0 & 0 & 1 & 1\\
1 & 0 & 0 & 0 & 1\\
1 & 1 & 0 & 0 & 0
\end{array}
\right]
\]
satisfies $s(A)=3.$ If $A\in\mathcal{PC}_{6}^{(1)}$, then $s(A)$ can be
$0,1,2,3$ or $4$. If $A\in\mathcal{PC}_{7}^{(1)}$, then $s(A)$ can be any
integer from $0$ to $7.$ For example, any matrix $A\in\mathcal{PC}_{7}^{(1)}$
such that
\[
T(A)=\left[
\begin{array}
[c]{ccccccc}%
0 & 1 & 0 & 1 & 1 & 0 & 0\\
0 & 0 & 1 & 1 & 0 & 1 & 0\\
1 & 0 & 0 & 1 & 1 & 0 & 0\\
0 & 0 & 0 & 0 & 1 & 1 & 1\\
0 & 1 & 0 & 0 & 0 & 1 & 1\\
1 & 0 & 1 & 0 & 0 & 0 & 1\\
1 & 1 & 1 & 0 & 0 & 0 & 0
\end{array}
\right]
\]
satisfies $s(A)=7.$

\bigskip

\begin{proposition}
\label{PuniqT}Let $A\in\mathcal{PC}_{n}^{(1)}$ and suppose that $s(A)=1$. Let
$B=[b_{ij}]$ and $B^{\prime}=[b_{ij}^{\prime}]$ be s-minimizers of $A$, and
$T(B)$ and $T(B^{\prime})$ be the corresponding tournament matrices. Then
$T(B)=T(B^{\prime})$.
\end{proposition}

\begin{proof}
Suppose that $s(A)=1$. Let $j>i$ and $l>k$ be such that $b_{ij}<1$ and
$b_{kl}^{\prime}<1.$ Then, all entries above the diagonal of $T(B)$ and
$T^{\prime}(B)$ are $1$, except those in positions $i,j$ and $k,l$,
respectively$,$ which are $0$. Also, all entries bellow the diagonal of $T(B)$
and $T^{\prime}(B)$ are $0$, except those in positions $j,i$ and $l,k$,
respectively, which are $1$. Notice that, by Lemma \ref{PP1}, $j\geq i+2$ and
$l\geq k+2$.

Since row $q$ of an $n$-by-$n$ matrix has $n-q$ entries above the diagonal,
and the diagonal entries of $T(B)$ are zero, the row sums corresponding to
rows $i$ and $i+1$ of $T(B)$ are both $n-i-1$, the row sums corresponding to
rows $j$ and $j-1$ of $T(B)$ are both $n-j+1$, and the sum of row $p$, with
$p\neq i,i+1,j-1,j$, is $n-p.$ Thus, the set of row sums in $T(B)$ is exactly
$\{0,\ldots,n-1\}\backslash\{n-i,n-j\}$.

Similarly, the set of row sums in $T(B^{\prime})$ is exactly $\{0,\ldots
,n-1\}\backslash\{n-k,n-l\}$.

As $B$ and $B^{\prime}$ are s-minimizers of $A$, they are permutationally
similar (as both are permutationally similar to $A$). Thus, they have the same
set of row sums. By hypothesis, $n-i>n-j$ and $n-k>n-l$, which then implies
$i$ must be $k$ and $j$ must be $l$, for the row sum sets to be as stated.
This means that $T(B)=T(B^{\prime})$.
\end{proof}

\bigskip

If $s(A)>1$ (which implies $n>4$), we may have two permutationally similar
s-minimizers $B$ and $B^{\prime}$ with $T(B)\neq T(B^{\prime}).$

\begin{example}
\label{extwomin}Let
\[
B=\left[
\begin{array}
[c]{ccccc}%
1 & 2 & \frac{1}{2} & 2 & 2\\
\frac{1}{2} & 1 & 2 & 2 & 2\\
2 & \frac{1}{2} & 1 & 2 & \frac{1}{2}\\
\frac{1}{2} & \frac{1}{2} & \frac{1}{2} & 1 & 2\\
\frac{1}{2} & \frac{1}{2} & 2 & \frac{1}{2} & 1
\end{array}
\right]  \text{ and }B^{\prime}=\left[
\begin{array}
[c]{ccccc}%
1 & 2 & \frac{1}{2} & 2 & \frac{1}{2}\\
\frac{1}{2} & 1 & 2 & 2 & 2\\
2 & \frac{1}{2} & 1 & 2 & 2\\
\frac{1}{2} & \frac{1}{2} & \frac{1}{2} & 1 & 2\\
2 & \frac{1}{2} & \frac{1}{2} & \frac{1}{2} & 1
\end{array}
\right]  ,
\]
both in $\mathcal{PC}_{5}^{(1)}$. We have $B^{\prime}=P^{T}BP$, with%
\[
P=\left[
\begin{array}
[c]{ccccc}%
0 & 1 & 0 & 0 & 0\\
0 & 0 & 1 & 0 & 0\\
1 & 0 & 0 & 0 & 0\\
0 & 0 & 0 & 1 & 0\\
0 & 0 & 0 & 0 & 1
\end{array}
\right]  .
\]
Now, $s(B)=2$. It is at most $2$ because $T(B)$ has $2$ zeros above the
diagonal. And it is at least $2$ because the tournament has $2$ independent
3-cycles, which ensures that any permutation similarity has at least $2$ zeros
above the diagonal. Thus, $B$ and $B^{\prime}$ are both s-minimizers and
$T(B)\neq T(B^{\prime})$.

For any size $n>5$, we may construct two permutationally similar s-minimizers
$C$ and $C^{\prime}$, with $s(C)=2$ and $T(C)\neq T(C^{\prime})$. Simply add
$n-5$ rows on the top of $B$ and $B^{\prime}$, and the corresponding columns
on the left, so that the resulting reciprocal matrices $C$ and $C^{\prime}$
have all entries above the diagonal in the first $n-5$ rows $>$ $1$.
\end{example}

Even when the tournament matrix associated with an s-minimizer of $A$ is
unique, there may be different permutation similarities of $A$ that give an s-minimizer.

\begin{example}
\label{exsmin2}Let $A\in\mathcal{PC}_{3}^{(1)}$ with corresponding tournament
matrix%
\begin{equation}
T(A)=\left[
\begin{array}
[c]{ccc}%
0 & 1 & 0\\
0 & 0 & 1\\
1 & 0 & 0
\end{array}
\right]  . \label{exTA}%
\end{equation}
Since $s(A)$ cannot be $0$, we have $s(A)=1$. Consider the permutation matrix
\begin{equation}
P=\left[
\begin{array}
[c]{ccc}%
0 & 1 & 0\\
0 & 0 & 1\\
1 & 0 & 0
\end{array}
\right]  . \label{exTA2}%
\end{equation}
We have $T(A)=T(PAP^{T})=T(P^{T}AP)$. Thus, $A$, $PAP^{T}$ and $P^{T}AP$ are
s-minimizers. Moreover, no other permutation similarity of $A$ is an s-minimizer.

For example, if
\[
A=\left[
\begin{array}
[c]{ccc}%
1 & 2 & \frac{1}{2}\\
\frac{1}{2} & 1 & 3\\
2 & \frac{1}{3} & 1
\end{array}
\right]  \in\mathcal{PC}_{3}^{(1)},
\]
then%
\begin{equation}
B:=PAP^{T}=\left[
\begin{array}
[c]{ccc}%
1 & 3 & \frac{1}{2}\\
\frac{1}{3} & 1 & 2\\
2 & \frac{1}{2} & 1
\end{array}
\right]  \label{AA2}%
\end{equation}
is another s-minimizer, though $T(A)=T(B)$. If
\begin{equation}
A=\left[
\begin{array}
[c]{ccc}%
1 & 2 & \frac{1}{2}\\
\frac{1}{2} & 1 & 2\\
2 & \frac{1}{2} & 1
\end{array}
\right]  \in\mathcal{PC}_{3}^{(1)}, \label{A1}%
\end{equation}
we have $P^{T}AP=PAP^{T}=A.$ Thus, $A$ is the only s-minimizer.
\end{example}

\subsection{The best weight vectors\label{sbest}}

Linear programming methods have been offered to minimize the number of rank
reversals given by a weight vector relative to a reciprocal matrix, see
\cite{Fara2,Gass,TuWu,wang}, for example. Here we address this problem via a
direct theoretical approach. \bigskip

First we characterize those strict vectors $w\in\mathbb{R}_{+}^{n}$ (not
necessarily efficient) that give the fewest rank reversals. We then show the
existence of efficient such vectors. Notice that vectors with the same strict
order give the same rank reversals.

\begin{lemma}
\label{tminrank}Suppose that $B\in\mathcal{PC}_{n}^{(1)}$ is an s-minimizer.
Then, the fewest rank reversals for $B$ given by a strict $w\in\mathbb{R}%
_{+}^{n}$ is $2s(B)$ ($s(B)$ rank reversals above, and $s(B)$ below, the
diagonal; so $s(B)$ is the number of pairs of rank reversals). Moreover, if
$w\in\mathbb{R}_{+}^{n}$ is strictly decreasing, then $w$ gives $s(B)$ pairs
of rank reversals for $B$.
\end{lemma}

\begin{proof}
Clearly, the number of rank reversals for $B$ given by a strictly decreasing
vector is $2s(B)$. Suppose that there is a strict vector $w$ that gives
$2s^{\prime}$, with $s^{\prime}<s(B)$, rank reversals for $B$, and let $P$ be
a permutation matrix such that $Pw$ is strictly decreasing. By Proposition
\ref{Prankpermut}, $Pw$ gives $2s^{\prime}$ rank reversals for $PBP^{T}$.
Thus, $PBP^{T}$ has $s^{\prime}$ entries $<1$ above the diagonal, since all
entries above the diagonal of $(Pw)(Pw)^{(-T)}$ are $>1$. This gives a
contradiction because $B$ is an s-minimizer and $s^{\prime}<s(B)$.
\end{proof}

\bigskip The rank reversals for $B$, given by a strictly decreasing vector,
occur in the positions above the diagonal with entries $<1$, and in the
corresponding positions below the diagonal.

\bigskip

Without our prevailing assumptions of $\mathcal{PC}_{n}^{(1)}$ and strict
vectors, the tally of rank reversal would be different. For example, the
vector of ones would have no rank reversals.

\bigskip

We have the following consequence of Proposition \ref{Prankpermut} and Lemma
\ref{tminrank}, together with a simple observation. If a strictly decreasing
vector gives $s(B)$ pairs of rank reversals for $B\in\mathcal{PC}_{n}^{(1)}$,
then $B$ is an s-minimizer.

\begin{theorem}
\label{CC1}Let $A\in\mathcal{PC}_{n}^{(1)}$. The fewest rank reversals for $A$
given by a strict $w\in\mathbb{R}_{+}^{n}$ is $2s(A)$. Moreover, if
$w\in\mathbb{R}_{+}^{n}$ is a strict vector, then $w$ has $s(A)$ pairs of rank
reversals for $A$ if and only if $PAP^{T}$ is an s-minimizer, in which $P$ is
the $n$-by-$n$ permutation matrix such that $Pw$ is decreasing.
\end{theorem}

We can now describe the number of strict orders that give the fewest rank reversals.

\begin{corollary}
\label{CC2}Let $A\in\mathcal{PC}_{n}^{(1)}$. The number of orders among the
strict vectors in $\mathbb{R}_{+}^{n}$ that give $s(A)$ pairs of rank
reversals for $A$ is the number of $n$-by-$n$ permutation matrices $P$ such
that $PAP^{T}$ is an s-minimizer.
\end{corollary}

The notion of a rank reversal for $A\in\mathcal{PC}_{n}^{(1)}$ can be given in
terms of the tournament matrix $T(A)=[t_{ij}]$ associated with $A.$ A strict
$w\in\mathbb{R}_{+}^{n}$ has an $i,j$ rank reversal for $A$ if and only if
$t_{ij}=1$ and $w_{j}>w_{i}$ ($j$ is ranked higher than $i$) or $t_{ij}=0$ and
$w_{i}>w_{j}$ ($i$ is ranked higher than $j$). With this definition in mind,
Lemma \ref{tminrank}, Theorem \ref{CC1} and Corollary \ref{CC2} could have
been discussed in the context of tournaments. See, for example, \cite{Ali} for
a discussion of related ideas.

\bigskip

Now we use Theorem \ref{CC1} to obtain our result for efficient vectors. It
implies that, for any order of a weight vector that gives the fewest rank
reversals for $A$, there always exists an efficient vector for $A$ with that
order. So, Corollary \ref{CC2} also applies to efficient vectors.

\begin{theorem}
\label{thrr}Let $A\in\mathcal{PC}_{n}^{(1)}$. There is a strict vector in
$\mathcal{E}(A)$ with $s(A)$ pairs of rank reversals, that is, with the fewest
rank reversals among all strict vectors in $\mathbb{R}_{+}^{n}$. If $A$ is an
s-minimizer, there is a decreasing such vector.
\end{theorem}

\begin{proof}
Let $P$ be a permutation matrix such that $B=PAP^{T}\ $is an $s$-minimizer.
Since, by Lemma \ref{PP1}, the entries on the superdiagonal of $B$ are $>1$,
by Proposition \ref{cdecrsuper}, there is a strictly decreasing vector
$w\in\mathcal{E}(B)$. Thus, by Theorem \ref{CC1}, $P^{T}w$ has $s(A)$ pairs of
rank reversals for $A$. By Lemma \ref{lsim}, $P^{T}w\in\mathcal{E}(A)$.
\end{proof}

\bigskip

The strict efficient vector with $s(A)$ pairs of rank reversals, guaranteed by
the proof of Theorem \ref{thrr}, is obtained from a certain Hamiltonian path.
Sufficiently small perturbations of this vector that are in $\mathcal{E}(A)$
are also strict and have $s(A)$ pairs of rank reversals.

\begin{example}
Consider the matrix $A\in\mathcal{PC}_{4}^{(1)}$ in Example \ref{Exnoamb}.
Then,
\[
B:=PAP^{T}=\left[
\begin{array}
[c]{cccc}%
1 & 2 & 2 & \frac{1}{2}\\
\frac{1}{2} & 1 & 8 & 2\\
\frac{1}{2} & \frac{1}{8} & 1 & 2\\
2 & \frac{1}{2} & \frac{1}{2} & 1
\end{array}
\right]  \text{, with }P=\left[
\begin{array}
[c]{cccc}%
0 & 1 & 0 & 0\\
1 & 0 & 0 & 0\\
0 & 0 & 0 & 1\\
0 & 0 & 1 & 0
\end{array}
\right]  .
\]
Thus, $s(A)=s(B)=1,$ as $B$ has no row with all off-diagonal entries $>1$ and
has just one entry $<1$ above the diagonal. According to the proof of
Proposition \ref{cdecrsuper}, let $w$ be such that $ww^{(-T)}$ coincides with
$B$ along the superdiagonal, that is, $\frac{w_{1}}{w_{2}}=\frac{w_{3}}{w_{4}%
}=2$ and $\frac{w_{2}}{w_{3}}=8.$ Projectively, this vector is $w=\left[
\begin{array}
[c]{cccc}%
32 & 16 & 2 & 1
\end{array}
\right]  ^{T}.$ Then, $w\in\mathcal{E}(B)$. By Lemma \ref{tminrank} (or check
directly), $w$ has the fewest pairs of rank reversals for $B$, which is
$s(B)=1.$ Then, $P^{T}w=\left[
\begin{array}
[c]{cccc}%
16 & 32 & 1 & 2
\end{array}
\right]  ^{T}$ gives the fewest rank reversals for $A$. Also, by Lemma
\ref{lsim}, $P^{T}w\in\mathcal{E}(A)$.
\end{example}

If there are two different s-minimizers for the same matrix (see Example
\ref{extwomin}, or Example \ref{exsmin2} with $A$ being the matrix
(\ref{AA2})), they are permutationally similar. In this case, because of
Corollary \ref{CC2}, there exists more than one strict efficient order with
$s(A)$ pairs of rank reversals. This also occurs if there is just one
s-minimizer, but the permutation similarity that produces it is not unique.

\begin{example}
Let $A\in\mathcal{PC}_{3}^{(1)}$, with $T(A)$ as in (\ref{exTA}), and let $P$
be as in (\ref{exTA2}). Taking into account Theorem \ref{thrr} and the
observations in Example \ref{exsmin2}, it follows that there are strictly
decreasing vectors $z\in\mathcal{E}(A)$, $u\in\mathcal{E}(PAP^{T})$ and
$v\in\mathcal{E}(P^{T}AP)$, each of them giving one pair of rank reversals for
the corresponding matrix. By Lemma \ref{lsim}, $P^{T}u,Pv\in\mathcal{E}(A)$,
and, by Proposition \ref{Prankpermut}, each of these vectors give one pair of
rank reversals for $A$. Thus, there is an efficient vector $w$ with strict
decreasing order and giving one pair of rank reversals, one such that
$w_{2}>w_{3}>w_{1},$ and one such that $w_{3}>w_{1}>w_{2}$. Thus, the order of
the vectors in $\mathcal{E}(A)$ with $s(A)$ pairs of rank reversals is not
unique, even if the s-minimizer is unique (for example if $A$ is the matrix
(\ref{A1})).
\end{example}

A \emph{partial reciprocal matrix} is a matrix with unspecified entries,
occurring in symmetrically placed positions, and that allows a completion to a
reciprocal matrix. Given a partial reciprocal matrix $A$, we may want to
complete it to a reciprocal matrix $A^{\prime}$ with minimum $s(A^{\prime})$.
Such a completion is obtained as follows. We consider a matrix permutationally
similar to the partial matrix $A$ with the fewest specified entries $<1$ above
the diagonal. Call this number $s$. Now complete the unspecified entries above
the diagonal of this permutation similarity of $A$ with numbers $>1$ and the
corresponding entries below the diagonal in such a way that a reciprocal
matrix $B$ is obtained. Then $s(B)=s$. By Lemma \ref{tminrank}, any strictly
decreasing vector $w$, which can be chosen efficient by Theorem \ref{thrr},
gives $s(B)$ pairs of rank reversals for $B$. Therefore, the permutation of
$w$ described in Theorem \ref{CC1} gives $s(B)$ pairs of rank reversals for
$A^{\prime}$. By Lemma \ref{lsim}, this vector is efficient for $A^{\prime}$
if $w$ is efficient for $B$. Moreover, it has no rank reversals for the
completion $A^{\prime}$ in the positions of the unspecified entries in $A$ and
has the minimum possible number of rank reversals for the specified entries.

\section{Order transitive reciprocal matrices\label{stransitive}}

We call $A\in\mathcal{PC}_{n}^{(1)}$ \emph{order transitive} if $a_{ij}>1$ and
$a_{jk}>1$ implies $a_{ik}>1\ $for all $1\leq i,j,k\leq n$. If $A\ $is order
transitive, we also have that $a_{ij}<1$ and $a_{jk}<1$ implies $a_{ik}<1.$

\begin{example}
Let
\[
A=\left[
\begin{array}
[c]{ccc}%
1 & 2 & 3\\
\frac{1}{2} & 1 & 4\\
\frac{1}{3} & \frac{1}{4} & 1
\end{array}
\right]  \text{ and }B=\left[
\begin{array}
[c]{ccc}%
1 & 2 & \frac{1}{2}\\
\frac{1}{2} & 1 & 4\\
2 & \frac{1}{4} & 1
\end{array}
\right]  ,
\]
both in $\mathcal{PC}_{3}^{(1)}$. Matrix $A$ is order transitive. Matrix $B$
is not, as $a_{31}=2$, $a_{12}=2$ and $a_{32}=\frac{1}{4}<1.$
\end{example}

Order transitivity is invariant under permutation similarity.

\bigskip

The outcomes of the tournament matrix $T(A)$ associated with an order
transitive $A\in\mathcal{PC}_{n}^{(1)}$ are transitive and each contest result
may be deduced from the rank order of the contestants. With this connection,
the following proposition follows from the theory of tournaments \cite{Moon}.

\begin{proposition}
\label{p1}Let $A\in\mathcal{PC}_{n}^{(1)}.$ Matrix $A$ is order transitive if
and only if it is permutationally similar to a matrix with all entries above
the diagonal $>1$ (that is, if and only if $s(A)=0$).
\end{proposition}

\bigskip

A statement corresponding to the "if" claim in Proposition \ref{p1} does not
hold if $A\in\mathcal{PC}_{n}^{0}$ has off-diagonal entries equal to $1$ and
the condition defining order transitivity is changed to $a_{ij}\geq1$ and
$a_{jk}\geq1$ implies $a_{ik}\geq1$ (as in \cite{Cavallo,TuWu}).

\begin{example}
Let%
\[
C=\left[
\begin{array}
[c]{ccc}%
1 & 2 & 1\\
\frac{1}{2} & 1 & 4\\
1 & \frac{1}{4} & 1
\end{array}
\right]  \in\mathcal{PC}_{3}^{0}.
\]
Matrix $C=[c_{ij}]$ has all upper diagonal entries $\geq1$. However,
$c_{31}\geq1$, $c_{12}\geq1$ and $c_{32}<1$.
\end{example}

An order transitive matrix $A\in\mathcal{PC}_{n}^{(1)}$ is easy to identify.
It is such a matrix if and only if, for any $i\in\{1,\ldots,n-1\}$, there is a
row with exactly $i$ entries $>1$. Any matrix $A\in\mathcal{PC}_{n}^{(1)}$
whose columns have the same order (see Section \ref{s81}) is order transitive.
This is because $A$ is permutationally similar to a matrix $B$ with decreasing
columns. Since the diagonal entries are $1,$ all entries of $B$ above the
diagonal are $>1$.\bigskip

For a general $A\in\mathcal{PC}_{n}^{(1)}$, the number $s(A)$ is the smallest
number of pairs of symmetrically placed entries in $A$ that need to be changed
in order to get an order transitive matrix. More precisely, if $B$ is an
s-minimizer of $A$, an appropriate modification of the pairs of symmetrically
placed entries in $B$ whose entry above the diagonal is $<1$ (for example the
interchange of the entries in the pair) gives an order transitive matrix.
Reversal of the permutation that produces $B$, when applied to $B$, gives an
order transitive matrix in place of $A$.

\bigskip

It is known (and also follows from the results here) that $A\in\mathcal{PC}%
_{n}^{(1)}$ is order transitive if and only if there is a strict vector in
$\mathbb{R}_{+}^{n}$ with no rank reversals \cite{Cavallo}. However, no
mention of the existence of an efficient vector was made. This is proven in
the next corollary. Notice that, if $A$ is order transitive, there is a unique
order for a strict vector with no rank reversals for $A$. This follows from
Corollary \ref{CC2}, since, if $A$ is order transitive, there is a unique
s-minimizer and it is attained by a unique permutation similarity.

\begin{corollary}
\label{ctrans}Let $A\in\mathcal{PC}_{n}^{(1)}$. There is a strict
$w\in\mathcal{E}(A)$ with no rank reversals for $A$ if and only if $A$ is
order transitive.
\end{corollary}

\begin{proof}
If there is a strict $w\in\mathcal{E}(A)$ with no rank reversals for $A$, by
Theorem \ref{CC1}, $s(A)=0$. By Proposition \ref{p1}, $A$ is order transitive.

The converse claim follows from Theorem \ref{thrr}, as, if $A$ is order
transitive, $s(A)=0$.
\end{proof}

\bigskip

We make additional observations about rank reversals in an order transitive matrix.

\bigskip

If $A\in\mathcal{PC}_{n}^{(1)}$ is order transitive, it has an entry
$a_{ij}>1$ along any cycle $\gamma\in\Gamma(A)$. Thus, because $P_{A,\gamma
}(i,j)=a_{ij}>1$, by Theorems \ref{Porderambgamma} and
\ref{Porderreversalgamma}, there is no $i,j$ rank reversal in $\mathcal{E}%
_{\gamma}(A)$.

\begin{theorem}
\label{torrr}Let $A\in\mathcal{PC}_{n}^{(1)}$ be order transitive and let
$\gamma$ be the only H-cycle in $A$ with $n-1$ entries $>1$. If $\gamma
(A)\in\Gamma(A)$, then no vector in $\mathcal{E}_{\gamma}(A)$ has a rank
reversal in any position.
\end{theorem}

\begin{proof}
Because of Propositions \ref{Prankpermut} and \ref{p1}, we assume that $A$ has
all entries above the diagonal $>1$. Then $\gamma=12\cdots n1.$ If
$\gamma(A)\in\Gamma(A)$, by Corollary \ref{ctr}, all vectors in $\mathcal{E}%
_{\gamma}(A)$ are strictly decreasing. Because all entries above the diagonal
of $A$ are $>1$, no vector in $\mathcal{E}_{\gamma}(A)$ has a rank reversal.
\end{proof}

\begin{example}
\label{exUO}Continue Examples \ref{exbound0}, \ref{exbound02},
\ref{Exnoorderamb}, \ref{exbound03} and \ref{exrr}. Matrix $A$ is order
transitive. All vectors in $\mathcal{E}_{\alpha}(A)$ are strictly decreasing
and, thus, give no rank reversals for $A$. This is in agreement with Theorem
\ref{torrr}, as there are $n-1$ entries $>1$ in $A$ along the H-cycle
$\alpha\in\Gamma(A)$. The strict vector $w\in\mathcal{E}(A)$ that gives no
rank reversals for $A$, guaranteed by Theorem \ref{thrr} and constructed as in
the proof of Proposition \ref{cdecrsuper}, is defined by $\frac{w_{i}}%
{w_{i+1}}=a_{i,i+1}$, $i=1,2,3$, that is, $w=\left[
\begin{array}
[c]{cccc}%
8 & 4 & 2 & 1
\end{array}
\right]  ^{T}$ (projectively). We have $w\in\mathcal{E}_{\alpha}(A)$.
\end{example}

Rank reversals of the "same type" are "transitive".

\begin{proposition}
Let $A\in\mathcal{PC}_{n}^{(1)}$ be order transitive and $w\in\mathbb{R}%
_{+}^{n}$ be strict. Let $i,j,k\in\{1,\ldots,n\}$. If $a_{ij}>1$ and
$a_{jk}>1$ (resp. $a_{ij}<1$ and $a_{jk}<1$) and $w$ has a rank reversal at
positions $i,j$ and $j,k$, then $w$ has a rank reversal at position $i,k.$
\end{proposition}

\begin{proof}
If $w$ has a rank reversal at $i,j$ and $j,k$, and $a_{ij}>1$ and $a_{jk}>1$
(which implies that $i,j,k$ are distinct), then $\frac{w_{i}}{w_{j}}<1$ and
$\frac{w_{j}}{w_{k}}<1,$ implying $\frac{w_{i}}{w_{k}}<1$. Since $A$ is order
transitive, we have $a_{ik}>1$. Thus, $w$ has a rank reversal at position
$i,k$. The parenthetical claim follows similarly.
\end{proof}

\bigskip

Any number $s(A)$ possible for an $A\in\mathcal{PC}_{n}^{(1)}$ can be attained
for some $B\in\mathcal{PC}_{n}^{(1)}$ arbitrarily close to a consistent
matrix. Simply consider the entrywise power of $A$ with sufficiently small
exponent $\varepsilon>0$. Thus, with respect to rank reversals and in light of
Theorem \ref{CC1}, this is in agreement with the observation in \cite{TuWu}
that, more important than a low index of inconsistency of a reciprocal matrix,
is how close the matrix is to order transitivity.

\section{Uniformly ordered reciprocal matrices\label{sUO}}

In Section \ref{stransitive} we have shown that $A\in\mathcal{PC}_{n}^{(1)}$
is order transitive if and only if there is an efficient vector for $A$ with
no rank reversals. Here we study the class of reciprocal matrices in
$\mathcal{PC}_{n}^{(1)}$ for which no efficient vector has a rank reversal.
This is contained in the family of order transitive matrices. We also study
the case in which all vectors in $\mathcal{E}_{\gamma}(A)$, for some
$\gamma\in\Gamma(A)$, have the same rank reversals relative to $A$. Most of
the results can be stated for matrices in $\mathcal{PC}_{n}^{0}$.

\subsection{Definitions and basic properties\label{s81}}

Because we have heretofore focused on strict vectors, we first clarify what is
meant by "order" for vectors in $\mathbb{R}_{+}^{n}$ with possible equal entries.

We say that vectors $v^{(1)},\ldots,v^{(k)}\in\mathbb{R}_{+}^{n}$ have the
\emph{same order} if there is an $n$-by-$n$ permutation matrix $Q$ such that
$Qv^{(i)}$ is (weakly) decreasing for all $i=1,\ldots,k$. If one of the
vectors is strict, the matrix $Q$ is unique, in which case there is just one
(weak) order shared by the vectors.

\begin{example}
The vectors
\[
\left[
\begin{array}
[c]{ccc}%
3 & 2 & 1
\end{array}
\right]  ^{T}\text{ and }\left[
\begin{array}
[c]{ccc}%
4 & 1 & 1
\end{array}
\right]  ^{T}%
\]
have the same order. The common order is unique (decreasing). The vectors
\[
\left[
\begin{array}
[c]{ccc}%
3 & 2 & 2
\end{array}
\right]  ^{T}\text{ and }\left[
\begin{array}
[c]{ccc}%
4 & 1 & 1
\end{array}
\right]  ^{T}%
\]
have the same order. The common order is not unique. It may be the weak
decreasing order or the (weak) order first, third, second.
\end{example}

Matrix $A\in\mathcal{PC}_{n}^{0}$ is called \emph{uniformly ordered} (UO) if
all the vectors in $\mathcal{E}(A)$ have the same order, and uniformly ordered
for $\gamma\in\Gamma(A)$, if all the vectors in $\mathcal{E}_{\gamma}(A)$ have
the same order. The former definition was introduced in \cite{FJ20} and the
latter in \cite{FJ8}. Of course, if $A$ is UO, then it is UO for any
$\gamma\in\Gamma(A).$

\bigskip

Since any column of $A\in\mathcal{PC}_{n}^{0}$ is efficient for $A$
\cite{FJ1}, if $A$ is UO, all columns have the same order. Thus, a UO matrix
$A\in\mathcal{PC}_{n}^{(1)}$ is order transitive.

\bigskip

We show that there is only one order shared by all columns of a UO matrix. To
see this, we first show that such a matrix cannot have equal rows and, thus,
equal columns. Of course, this is clear if the matrix is in $\mathcal{PC}%
_{n}^{(1)}$, but interestingly it is still true in $\mathcal{PC}_{n}^{0}$. The
following lemmas are helpful.

\begin{lemma}
\label{leqrows}Let $V$ be the matrix whose columns are $v^{(1)},\ldots
,v^{(k)}\in\mathbb{R}_{+}^{n}$, which we assume have the same (weak) order.
The permutation matrix $Q$ such that all the columns of $QV$ are (weakly)
decreasing is unique if and only if $V$ has no equal rows.
\end{lemma}

\begin{proof}
Suppose that $V$ has some equal rows. Without loss of generality, we may
assume that the last two rows are equal. Let $Q$ be such that all the columns
of $QV$ are (weakly) decreasing. Let
\[
P=I_{n-2}\oplus\left[
\begin{array}
[c]{cc}%
0 & 1\\
1 & 0
\end{array}
\right]  .
\]
Then, $Q\neq QP$ and all columns of $QV$ and $QPV$ are decreasing.

To prove the converse, suppose that $V$ has no equal rows. Moreover, suppose,
without loss of generality, that the columns of $V$ are decreasing. Let $Q\neq
I_{n}$ be a permutation matrix. Let $i\in\{1,\ldots,n-1\}$ be the smallest
index for which there is $j\neq i$, such that the $i,j$ entry in $Q$ is $1$
(that is, $i$ is the smallest index of a row of $Q$ with a $1$ in an
off-diagonal position). Note that $j>i$. Let $r_{i}$ and $r_{j}$ denote the
$i$th and $j$th rows of $V$, respectively. By assumption $r_{i}\neq r_{j}%
\ $and the submatrix $\left[
\begin{array}
[c]{c}%
r_{i}\\
r_{j}%
\end{array}
\right]  $ of $V$ has decreasing columns, implying that there is at least one
strictly decreasing column. Thus, not all columns of the submatrix $\left[
\begin{array}
[c]{c}%
r_{j}\\
r_{i}%
\end{array}
\right]  $ of $QV$ are decreasing.
\end{proof}

\bigskip

The following lemma was given in \cite{FJblock}. We denote by $w[\{1,\ldots
,s\}]$ the vector formed by the first $s$ entries of $w\in\mathbb{R}_{+}^{n},$
$1\leq s\leq n.$ Recall that $J$ is a matrix of ones.

\begin{lemma}
\label{lcompl}Let
\begin{equation}
C=\left[
\begin{array}
[c]{cc}%
B & J\\
J & J
\end{array}
\right]  \in\mathcal{PC}_{n}^{0}, \label{BJ}%
\end{equation}
with $B\in\mathcal{PC}_{s}$, $s<n.$ Let $w=\left[
\begin{array}
[c]{cccc}%
w_{1} & w_{2} & \cdots & w_{n}%
\end{array}
\right]  ^{T}\in\mathbb{R}_{+}^{n},\ $with $w[\{1,\ldots,s\}]\in
\mathcal{E}(B).$ Then, $w\in\mathcal{E}(C)$ if and only if%
\[
\min_{1\leq i\leq s}w_{i}\leq w_{k}\leq\max_{1\leq i\leq s}w_{i}%
\qquad\text{for }s+1\leq k\leq n.
\]

\end{lemma}

\begin{theorem}
\label{Teqrows}Let $A\in\mathcal{PC}_{n}^{0}$ be UO. Then all rows (columns)
of $A$ are distinct.
\end{theorem}

\begin{proof}
Suppose that $A$ has two equal rows. Without loss of generality, suppose they
are the last two. Then, there is an $(n-2)$-by-$(n-2)$ positive diagonal
matrix $D$ such that $C=\left(  D^{-1}\oplus I_{2}\right)  A\left(  D\oplus
I_{2}\right)  $ has the form (\ref{BJ}) with $s=n-2$. Since $A$ is
inconsistent, we have $B\neq J.$ Thus, there is a $w\in\mathcal{E}(B)$ with
non-constant entries (for example a certain column of $B$). Then, $\min_{1\leq
i\leq n-2}w_{i}<\max_{1\leq i\leq n-2}w_{i}$. By Lemma \ref{lcompl}, there are
vectors $u,v\in\mathcal{E}(C)$ such that $u[\{1,\ldots,n-2\}]=v[\{1,\ldots
,n-2\}]=w,$ $u_{n-1}>u_{n}$ and $v_{n-1}<v_{n}$. By Lemma \ref{lsim}, $\left(
D\oplus I_{2}\right)  u,\left(  D\oplus I_{2}\right)  v\in\mathcal{E}(A).$
Vectors $\left(  D\oplus I_{2}\right)  u$ and $\left(  D\oplus I_{2}\right)
v$ do not have the same order, implying that $A$ is not UO.
\end{proof}

\bigskip

From Theorem \ref{Teqrows} and Lemma \ref{leqrows}, we have the following.

\begin{corollary}
\label{cdecr}If $A\in\mathcal{PC}_{n}^{0}$ is UO, there is a unique (weak)
order shared by all the columns of $A$.
\end{corollary}

\bigskip

In \cite{FJ20} we used, without comment, the fact that the efficient vectors
for a UO matrix with (weakly) decreasing columns are decreasing. This is
clarified using Corollary \ref{cdecr}.

\subsection{Characterization of UO and rank reversals\label{sUORR}}

Corollary \ref{Thintgea} (resp. Theorem \ref{Thintgamma}) and Theorem
\ref{rattain} allow us to characterize the UO (resp. UO for $\gamma\in
\Gamma(A)$) matrices $A\in\mathcal{PC}_{n}^{0}$. Part of the next result was
shown in \cite{FJ80}. We give a full proof here.

\begin{theorem}
\label{Tuo}Let $A\in\mathcal{PC}_{n}^{0}$. The following are equivalent:

\begin{enumerate}
\item $A$ is uniformly ordered;

\item for all $i,j\in\{1,\ldots,n\}$, $i\neq j$, $1\notin(L_{A}(i,j),\frac
{1}{L_{A}(j,i)})$, that is, $\mathcal{E}(A)$ is $i,j$ order unambiguous;

\item there is a permutation $i_{1}i_{2}\cdots i_{n}$ of $\{1,\ldots,n\}$ such
that $L_{A}(i_{s},i_{s+1})\geq1,$ for $s=1,\ldots,n-1$;

\item $L_{A}$ is permutationally similar to a matrix with all entries above
the diagonal $\geq1$;

\item $L_{A}$ has exactly $\frac{n^{2}-n}{2}$ off-diagonal entries $\geq1$;

\item for all $i,j\in\{1,\ldots,n\}$, $i\neq j$, exactly one of $L_{A}(i,j)$
or $L_{A}(j,i)$ is $\geq1$.
\end{enumerate}
\end{theorem}

\begin{proof}
(1.$\Leftrightarrow$2.) If $1\in(L_{A}(i,j),\frac{1}{L_{A}(j,i)})$, $i\neq j$,
then $L_{A}(i,j)<1$ and $\frac{1}{L_{A}(j,i)}>1.$ By Theorem \ref{rattain},
there is $w\in\mathcal{E}(A)\ $such that $\frac{w_{i}}{w_{j}}=L_{A}(i,j)<1$
and there is $v\in\mathcal{E}(A)\ $such that $\frac{v_{i}}{v_{j}}=\frac
{1}{L_{A}(j,i)}>1.$ This implies that $A$ is not UO. Thus, 1. implies 2. If 2.
holds, by Corollary \ref{Thintgea}, for all $i,j$, $i\neq j$, either
$\frac{w_{i}}{w_{j}}\geq1$ for all $w\in\mathcal{E}(A)$, or $\frac{w_{i}%
}{w_{j}}\leq1$ for all $w\in\mathcal{E}(A)$. Thus, statement 1. holds.

(2.$\Leftrightarrow$6.) Follows from Theorem \ref{Porderambgamma}. Notice
that, because of (\ref{inLA}), at least one of $L_{A}(i,j)$ or $L_{A}(j,i)$ is
$<1$.

(6.$\Rightarrow$3.) Suppose that 6. holds. Consider the directed graph with
$n$ vertices and an edge $(i,j)$ if and only if $L_{A}(i,j)\geq1$. Because of
the assumption, this graph is a tournament. By R\'{e}dei's Theorem, every
tournament contains a Hamiltonian path. So, 3. holds.

(3.$\Rightarrow$4.) Suppose that 3. holds. Assume, without loss of generality,
that $L_{A}(i,i+1)\geq1$, for $i=1,\ldots,n-1.$ Then, by Corollary
\ref{Thintgea}, $\frac{w_{i}}{w_{i+1}}\geq1$ for all $i=1,\ldots,n-1$ and all
$w\in\mathcal{E}(A)$. Thus, any $w\in\mathcal{E}(A)$ is a (weakly) decreasing
vector. By Theorem \ref{rattain}, $L_{A}(i,j)\geq1$, for all $j>i$. Thus, 4. holds.

(4.$\Rightarrow$5.$\Rightarrow$6.) Follows because of (\ref{inLA}).
\end{proof}

\bigskip

Not all order transitive matrices are UO.

\begin{example}
Continue Examples \ref{exbound0}, \ref{exbound02}, \ref{Exnoorderamb},
\ref{exbound03}, \ref{exrr} and \ref{exUO}. Matrix $A$ is order transitive but
is not UO (despite the fact that all columns are weakly decreasing), as
$L_{A}$ has $5<\frac{4^{2}-4}{2}$ off-diagonal entries $\geq1$. A\ decreasing
vector in $\mathcal{E}(A)$ was given in Example \ref{exUO}. A vector in
$\mathcal{E}(A)$ with a different order is, for example, $\left[
\begin{array}
[c]{cccc}%
8 & 2 & 4 & 1
\end{array}
\right]  ^{T}\in\mathcal{E}_{\beta}(A)$.
\end{example}

Using similar arguments, now applying Theorem \ref{Thintgamma} instead of
Corollary \ref{Thintgea}, and (\ref{Pijji}) instead of (\ref{inLA}), we get a
corresponding result for $\mathcal{E}_{\gamma}(A)$, part of which appeared in
\cite{FJ8}.

\begin{theorem}
\label{Tugama}Let $A\in\mathcal{PC}_{n}^{0}$ and $\gamma\in\Gamma(A).$ The
following are equivalent:

\begin{enumerate}
\item $A$ is uniformly ordered for $\gamma$;

\item for all $i,j\in\{1,\ldots,n\}$, $i\neq j$, $1\notin(P_{A,\gamma
}(i,j),\frac{1}{P_{A,\gamma}(j,i)})$, that is, $\mathcal{E}_{\gamma}(A)$ is
$i,j$ order unambiguous;

\item there is a permutation $i_{1}i_{2}\cdots i_{n}$ of $\{1,\ldots,n\}$ such
that $P_{A,\gamma}(i_{s},i_{s+1})\geq1,$ for $s=1,\ldots,n-1$;

\item $P_{A,\gamma}$ is permutationally similar to a matrix with all entries
above the diagonal $\geq1$;

\item $P_{A,\gamma}$ has exactly $\frac{n^{2}-n}{2}$ off-diagonal entries
$\geq1$;

\item for all $i,j\in\{1,\ldots,n\}$, $i\neq j$, exactly one of $P_{A,\gamma
}(i,j)$ or $P_{A,\gamma}(j,i)$ is $\geq1$.
\end{enumerate}
\end{theorem}

\bigskip

If $A$ is UO then, for each pair $i,j,$ $i\neq j$, the $i,j$ order unambiguity
in $\mathcal{E}_{\gamma}(A)$ is of the same type for all $\gamma\in\Gamma(A)$,
that is, $P_{A,\gamma}(i,j)\geq1$ for all $\gamma\in\Gamma(A),$ or
$P_{A,\gamma}(j,i)\geq1$ for all $\gamma\in\Gamma(A)$.$\ $Otherwise, there
would be an $i,j$ order ambiguity in $\mathcal{E}(A)$.

\bigskip

Matrix $A\in\mathcal{PC}_{n}^{0}$ is UO (resp. UO for $\gamma,$ with
$\gamma\in\Gamma(A)$), with all vectors in $\mathcal{E}(A)$ (resp.
$\mathcal{E}_{\gamma}(A)$) decreasing, if and only if $L_{A}(j,j+1)\geq1$
(resp. $P_{A,\gamma}(j,j+1)\geq1$) for all $j\in\{1,\ldots,n-1\}$.

\begin{example}
\label{Exdec}Let
\[
A=\left[
\begin{array}
[c]{cccc}%
1 & 2 & 8 & 32\\
\frac{1}{2} & 1 & 2 & 8\\
\frac{1}{8} & \frac{1}{2} & 1 & 2\\
\frac{1}{32} & \frac{1}{8} & \frac{1}{2} & 1
\end{array}
\right]  \in\mathcal{PC}_{4}^{(1)}.
\]
We have $\Gamma(A)=\{\gamma\}$, with $\gamma=12341$ and $\gamma(A)=\frac{1}%
{4}.$ Moreover,%
\[
L_{A}=P_{A,\gamma}=\left[
\begin{array}
[c]{cccc}%
1 & 2 & 4 & 8\\
\frac{1}{8} & 1 & 2 & 4\\
\frac{1}{16} & \frac{1}{8} & 1 & 2\\
\frac{1}{32} & \frac{1}{16} & \frac{1}{8} & 1
\end{array}
\right]  .
\]
By Theorem \ref{Tuo}, $A$ is UO since all entries on the superdiagonal of
$L_{A}$ are $\geq1$. Moreover, all efficient vectors for $A$ have decreasing
order. Since $L_{A}(i,i+1)>1$ for $i=1,2,3$, this order is strict. Since
$\#\Gamma(A)=1$, the fact that $A$ is UO also follows from Corollary \ref{ctr}.
\end{example}

See \cite[Example 26]{FJ8} for an example of a UO matrix $A$ with
$\#\Gamma(A)>1$.

\bigskip

If $\gamma\in\Gamma(A)$ is such that $A$ has $n-1$ entries $\geq1$ along
$\gamma$, then condition 3. in Theorem \ref{Tugama} holds and, thus, $A$ is UO
for $\gamma$. The converse is not generally true, as illustrated next. Note
that, because of R\'{e}dei's Theorem for semicomplete digraphs, for any
reciprocal matrix there is an H-cycle along which $A$ has at least $n-1$
entries $\geq1$, though this cycle is not necessarily in $\Gamma(A)$.

\begin{example}
Let
\[
A=\left[
\begin{array}
[c]{cccc}%
1 & \frac{1}{2} & 4 & 10\\
2 & 1 & 2 & 4\\
\frac{1}{4} & \frac{1}{2} & 1 & \frac{1}{2}\\
\frac{1}{10} & \frac{1}{4} & 2 & 1
\end{array}
\right]  \in\mathcal{PC}_{4}^{(1)}.
\]
We have $\gamma\in\Gamma(A)$, with $\gamma=13241$ and $\gamma(A)=\frac{4}{5}.$
Moreover,%
\[
P_{A,\gamma}=\left[
\begin{array}
[c]{cccc}%
1 & 2 & 4 & 8\\
\frac{2}{5} & 1 & \frac{8}{5} & 4\\
\frac{1}{5} & \frac{1}{2} & 1 & 2\\
\frac{1}{10} & \frac{1}{5} & \frac{2}{5} & 1
\end{array}
\right]  .
\]
By Theorem \ref{Tugama}, $A$ is UO for $\gamma$, though there are no $3$
entries $\geq1$ in $A$ along $\gamma.$ All vectors in $\mathcal{E}_{\gamma
}(A)$ are strictly decreasing.
\end{example}

\bigskip

From Theorems \ref{thallorderamb} and \ref{Tugama}, if $A\in\mathcal{PC}%
_{n}^{0}$ and there is an H-cycle $\gamma$ with all entries in $A$ along
$\gamma$ $<1$, then $A$ is not UO for $\gamma$, and, thus, is not UO. Notice
that $\gamma(A)<1,$ so that $\gamma\in\Gamma(A)$.

\bigskip

We have the following relationship between UO matrices and rank reversals.

\begin{corollary}
\label{cUOR}Let $A\in\mathcal{PC}_{n}^{(1)}.$ Then $A$ is UO if and only if no
vector in $\mathcal{E}(A)$ has a rank reversal in any position.
\end{corollary}

\begin{proof}
By Theorem \ref{Tuo}, $A$ is UO if and only if $\mathcal{E}(A)$ is $i,j$ order
unambiguous for any $i\neq j$. Now the result follows from Theorem
\ref{propanbrev}.
\end{proof}

\bigskip

Notice that, because of Corollary \ref{Thintgea} and Theorem \ref{rattain},
the "only if" claim in Corollary \ref{cUOR} holds if $A\in\mathcal{PC}_{n}%
^{0}$.

\section{Adding or deleting alternatives \label{sadddelete}}

Any $k$-by-$k$ principal submatrix $B$ of $A\in\mathcal{PC}_{n}^{0}$ lies in
$\mathcal{PC}_{k}$ and gives the pairwise comparisons among alternatives in
the subset corresponding to the rows (columns) of $B$. There is no simple
relationship between $\mathcal{E}(B)$ and $\mathcal{E}(A)$. However, cardinal
rankings given by an efficient vector for $B$ and an efficient vector for $A$
may happen to be the same for the alternatives in $B$. We explain how this
occurs. In that situation the rank reversals for the alternatives in $B$ given
by both vectors are the same.

\bigskip

Given a square matrix $A$ and a vector $w$, we denote by $A(i)$ the submatrix
of $A$ resulting from deletion of row and column $i$, and by $w(i)$ the vector
resulting from deletion of entry $i$.

\bigskip

\begin{theorem}
\label{text}\cite{FJ1} Let $A\in\mathcal{PC}_{n}^{0},$ $w=\left[
\begin{array}
[c]{cccc}%
w_{1} & w_{2} & \cdots & w_{n}%
\end{array}
\right]  ^{T}\in\mathbb{R}_{+}^{n},$ and $i\in\{1,\ldots,n\}$. Suppose that
\[
w(i)\in\mathcal{E}(A(i)).
\]
Then, $w\in\mathcal{E}(A)$ if and only if
\[
\min_{1\leq k\leq n,\text{ }k\neq i}\frac{w_{k}}{a_{ki}}\leq w_{i}\leq
\max_{1\leq k\leq n,\text{ }k\neq i}\frac{w_{k}}{a_{ki}}.
\]

\end{theorem}

For any principal submatrix $B$ of $A\in\mathcal{PC}_{n}^{0}$ and any
$v\in\mathcal{E}(B)$, application of Theorem \ref{text} guarantees that there
exist vectors $w\in\mathcal{E}(A)$ for which $v$ is a subvector of $w$ in the
"correct" positions (those corresponding to the rows of $B$). Thus, some
efficient vector for $A$ preserves the cardinal rankings in a given efficient
vector for $B$. Theorem \ref{text} was used to derive Lemma \ref{lcompl}.

\bigskip

Another, inductive, description of $\mathcal{E}(A)$ was given in \cite{FJ4}.
Let $\mathcal{E}(A;i)$ be the set of all vectors $w\in\mathcal{E}(A)$ such
that $w(i)\in\mathcal{E}(A(i))$.

\begin{theorem}
\label{thai}\cite{FJ4} Let $A\in\mathcal{PC}_{n}^{0}$, $n\geq4.$ Then,
\[
\mathcal{E}(A)=%
{\displaystyle\bigcup\limits_{i=1}^{n}}
\mathcal{E}(A;i)=%
{\displaystyle\bigcup\limits_{i=1,i\neq p}^{n}}
\mathcal{E}(A;i),
\]
for any $p\in\{1,\ldots,n\}$.
\end{theorem}

By Theorem \ref{thai}, if $w\in\mathcal{E}(A),$ there are at least two
distinct $i$'s such that $w(i)\in\mathcal{E}(A(i))$. For such $i$'s, the
cardinal rankings of the alternatives $1,\ldots,n$ without $i$ given by $w(i)$
are the same as those given by $w$. Successive application of the theorem
guarantees that, for any $3\leq s<n,$ there are subsets $S$ of $s$
alternatives such that the subvector of $w$ in the "correct" positions is
efficient for the principal submatrix $B$ of $A$ corresponding to the
alternatives in $S$. Thus, some efficient vector for $B$ preserves the
cardinal rankings of the alternatives in $S$ in the given efficient vector for
$A$.

\bigskip

Given $A\in\mathcal{PC}_{n}^{0}$, consider now the geometric convex hull of
the columns of $A$, $GC(A)$. In \cite{FJ2} it was shown that $GC(A)$ is
contained in $\mathcal{E}(A)$. Let $w\in GC(A)$ be obtained by taking
exponents equal to $0$ for a subset of columns indexed by a proper subset $S$
of $\{1,\ldots,n\}$. Then, $w\in\mathcal{E}(A)$. Let $B$ be the principal
submatrix resulting from $A$ by deleting the rows and columns indexed by $S$,
and let $v$ be obtained from $w$ by deleting the entries indexed by $S$. Then,
$v$ is a weighted geometric mean of the columns of $B$, with the same weights
as the nonzero weights used before. Then $v\in\mathcal{E}(B).$ For the
alternatives indexed by $\{1,\ldots,n\}\backslash S$, both vectors $w$ and $v$
give the same cardinal ranking.

\section{Numerical experiments\label{snum}}

Here we compare, in terms of number of rank reversals for a matrix
$A\in\mathcal{PC}_{n}^{0}$, the right Perron eigenvector (RP) \cite{saaty1977}%
$,$ the entrywise inverse of the left Perron eigenvector (LP)
\cite{CasatoLP,Johnson}, the geometric mean of RP and LP (GP) \cite{FJEJOR},
the geometric mean of the columns (GM) \cite{William}, a weighted geometric
mean of the columns with random weights (WGM) \cite{FJ2}, and a random
positive vector (RD). It is known that RP and LP may not be efficient
\cite{Boz,FJ5,FJEJOR}, while GP is universally efficient \cite{FJEJOR}. The
universal efficiency of WGM, shown in \cite{FJ2}, includes the case of
efficiency of GM, first shown in \cite{blanq2006}.

We generated $50000$ matrices $A$ in $\mathcal{PC}_{n}^{0}$. Each matrix $A$
is the Hadamard product of a random matrix $B$ with entries in $(0,1)$ and the
transpose of its entrywise inverse. As claimed in Section \ref{stransitive},
for our purpose, how close matrix $A$ is to consistency is not important, so
this method of generating $A$ works well.

We considered matrices of size $n$ from $2$ to $25.$ Since a rank reversal
occurs in position $i,j$ if and only if it occurs in position $j,i$, we just
counted the number of pairs with rank reversals. For each size $n$ and each
vector, we added the number of pairs with rank reversals in each trial and
divided this sum by the number of trials, so that the average of rank
reversals for each size was obtained. This number was then normalized by
division by the number of superdiagonal entries, $\frac{n^{2}-n}{2},$ giving
the average portion of pairs with rank reversals. In Table \ref{tab3} we give
the results obtained in percentage terms. It can be observed that GM is the
vector that performs the best, followed by WGM. Vector GP performs better than
RP and LP, which behave similarly. Vector RD has the worst behavior, as
expected. We also notice that the results get worse as $n$ increases.
\begin{table}[] \centering
\begin{tabular}
[c]{|c|c|c|c|c|c|c|}\hline
n & RP & LP & GP & GM & WGM & RD\\\hline
$2$ & $0$ & $0$ & $0$ & $0$ & $0$ & $49.62$\\\hline
$3$ & $22.36$ & $22.36$ & $22.36$ & $22.36$ & $21.43$ & $50.14$\\\hline
$4$ & $29.54$ & $29.63$ & $29.13$ & $27.62$ & $28.11$ & $49.95$\\\hline
$5$ & $34.02$ & $34.02$ & $33.26$ & $30.70$ & $31.68$ & $49.98$\\\hline
$6$ & $36.75$ & $36.86$ & $35.79$ & $32.67$ & $33.92$ & $50.02$\\\hline
$7$ & $38.88$ & $38.90$ & $37.76$ & $34.23$ & $35.48$ & $49.99$\\\hline
$8$ & $40.22$ & $40.27$ & $39.15$ & $35.25$ & $36.60$ & $49.96$\\\hline
$9$ & $41.40$ & $41.47$ & $40.31$ & $36.24$ & $37.52$ & $49.98$\\\hline
$10$ & $42.38$ & $42.39$ & $41.30$ & $37.05$ & $38.36$ & $50.01$\\\hline
$11$ & $43.04$ & $43.09$ & $42.01$ & $37.66$ & $38.94$ & $49.91$\\\hline
$12$ & $43.69$ & $43.69$ & $42.62$ & $38.20$ & $39.44$ & $49.95$\\\hline
$13$ & $44.18$ & $44.18$ & $43.16$ & $38.71$ & $39.93$ & $49.95$\\\hline
$14$ & $44.61$ & $44.59$ & $43.59$ & $39.14$ & $40.34$ & $50.04$\\\hline
$15$ & $44.99$ & $45.00$ & $44.01$ & $39.49$ & $40.67$ & $50.02$\\\hline
$16$ & $45.29$ & $45.28$ & $44.39$ & $39.86$ & $41.03$ & $49.99$\\\hline
$17$ & $45.61$ & $45.61$ & $44.71$ & $40.20$ & $41.33$ & $49.98$\\\hline
$18$ & $45.82$ & $45.85$ & $44.93$ & $40.45$ & $41.58$ & $49.98$\\\hline
$19$ & $46.07$ & $46.07$ & $45.22$ & $40.72$ & $41.81$ & $49.99$\\\hline
$20$ & $46.28$ & $46.26$ & $45.45$ & $40.97$ & $42.04$ & $49.99$\\\hline
$21$ & $46.45$ & $46.45$ & $45.66$ & $41.20$ & $42.25$ & $50.00$\\\hline
$22$ & $46.64$ & $46.64$ & $45.87$ & $41.43$ & $42.46$ & $50.03$\\\hline
$23$ & $46.78$ & $46.73$ & $46.00$ & $41.58$ & $42.59$ & $50.00$\\\hline
$24$ & $46.90$ & $46.90$ & $46.15$ & $41.76$ & $42.77$ & $50.00$\\\hline
$25$ & $47.00$ & $47.02$ & $46.30$ & $41.94$ & $42.90$ & $49.99$\\\hline
\end{tabular}
$\ \ \ \ \ \ \ \ \ \ \ \ \ \ $%
\caption{Average of the percentage of pairs with rank reversals  for 50000
random reciprocal matrices of size $n\leq25$, given by each of the vectors RP, LP, GP, GM, WGM and RD. }\label{tab3}%
\end{table}%

\section{Conclusions\label{concl}}

Our purpose was to better understand rank reversals between the cardinal
preference in the $i,j$ entry of a reciprocal matrix and those between the
same pair of alternatives in efficient vectors for the matrix, and to identify
those efficient vectors with the fewest rank reversals. We exploit recently
proven exact matrix bounds for the consistent matrices constructed from
efficient vectors. This allows us to characterize the existence or
nonexistence of a rank reversal in the $i,j$ entry for some efficient vector.
We determine all general weight vectors with the fewest rank reversals
relative to a given reciprocal matrix and how many reversals there are; and
demonstrate that there are always efficient vectors among these. Any order
that appears among the former also appears among the latter. These special
efficient vectors are candidates for the best weight vectors. The reciprocal
matrices for which there is an efficient vector with no rank reversals were
determined to be the "order transitive" matrices. The reciprocal matrices for
which no efficient vector has a rank reversal were studied in detail. These
matrices are called "uniformly ordered", as all efficient vectors have the
same order (so that there is no question of the best choice of order).
Additionally, the existence of efficient vectors with a specified order, for a
given reciprocal matrix, is studied. The preservation of cardinal rankings in
efficient vectors when adding, or deleting certain, alternatives is also
discussed. A number of other facts that may be useful are given. Numerical
experiments were conducted to study the performance in terms of rank reversals
of some often used formulaic weight vectors. We verified numerically that the
geometric mean, and the weighted geometric means, of the columns of a
reciprocal matrix, and the geometric mean of the right and entrywise inverse
left Perron eigenvector, behave better in terms of rank reversals than the
right Perron eigenvector and the entrywise inverse left Perron eigenvector,
individually. Interestingly, the former three vectors are universally
efficient while the latter are not.

\bigskip

\textbf{Declaration of conflict of interest}

\bigskip

There are no conflicts of interest.

\newpage


\bigskip

\begin{thebibliography}{99}                                                                                               %


\bibitem {Ali}I. Ali, W. D. Cook, M. Kress, \textit{On the minimum violations
ranking of a tournament}, Management Science 32 (1986), 660-672. https://doi.org/10.1287/mnsc.32.6.660.

\bibitem {Alon}N. Alon, \textit{Ranking Tournaments}, SIAM Journal on Discrete
Mathematics 20 (2006), 137-142. https://doi.org/10.1137/050623905.

\bibitem {Bana}C. A. Bana-e Costa, J. C. Vansnick, \textit{A critical analysis
of the eigenvalue method used to derive priorities in AHP}, European Journal
of Operational Research 187 (2008), 1422-1428. https://doi.org/10.1016/j.ejor.2006.09.022.

\bibitem {Barzilai-Cs}J. Barzilai, B. Golany, \textit{AHP rank reversal,
normalization and aggregation rules}, INFOR: Information Systems and
Operational Research 32 (1994), 57-64. https://doi.org/10.1080/03155986.1994.11732238.

\bibitem {BEL}V. Belton, T. Gear, \textit{On a short-coming of Saaty's method
of analytic hierarchies}, Omega 11 (1983), 228-230. https://doi.org/10.1016/0305-0483(83)90047-6.

\bibitem {blanq2006}R. Blanquero, E. Carrizosa, E. Conde, \textit{Inferring
efficient weights from pairwise comparison matrices}, Mathematical Methods of
Operations Research 64 (2006), 271-284. https://doi.org/10.1007/s00186-006-0077-1.

\bibitem {Boz}S. Boz\'{o}ki, \textit{Inefficient weights from pairwise
comparison matrices with arbitrarily small inconsistency},\textit{
}Optimization 63 (2014), 1893-1901. https://doi.org/10.1080/02331934.2014.903399.

\bibitem {Cavallo}B. Cavallo, L. D'Apuzzo, L. Basile, \textit{Weak consistency
for ensuring priority vectors reliability}, Journal of Multi-Criteria Decision
Analysis 23 (2016), 126-138. https://doi.org/10.1002/mcda.1559.

\bibitem {Chen}X. Chen, G. Ding, W. Zang, Q. Zhao, \textit{Ranking tournaments
with no errors I: Structural description}, Journal of Combinatorial Theory,
Series B 141 (2020), 264-294. https://doi.org/10.1016/j.jctb.2019.08.004.

\bibitem {Choo}E. U. Choo, W. C. Wedley, \textit{A common framework for
deriving preference values from pairwise comparison matrices}, Computers \&
Operations Research 31 (2004), 893-908. https://doi.org/10.1016/S0305-0548(03)00042-X.

\bibitem {William}G. Crawford, C. Williams, \textit{A note on the analysis of
subjective judgment matrices}, Journal of Mathematical Psychology 29 (1985),
387-405. https://doi.org/10.1016/0022-2496(85)90002-1.

\bibitem {Csato2017}L. Csat\'{o}, \textit{Eigenvector Method and Rank Reversal
in Group Decision Making Revisited}, Fundamenta Informaticae 156 (2017),
169-178. https://doi.org/10.3233/FI-2017-1602.

\bibitem {CasatoLP}L. Csat\'{o}, \textit{Right-left asymmetry of the
eigenvector method: A simulation study}, European Journal of Operational
Research 313 (2024), 708-717. https://doi.org/10.1016/j.ejor.2023.08.026.

\bibitem {Csato}L. Csat\'{o}, \textit{How to choose a completion method for
pairwise comparison matrices with missing entries: An axiomatic result},
International Journal of Approximate Reasoning 164 (2024), 109063. https://doi.org/10.1016/j.ijar.2023.109063.

\bibitem {Fara2}L. Faramondi, G. Oliva, S. Boz\'{o}ki, \textit{Incomplete
analytic hierarchy process with minimum weighted ordinal violations},
International Journal of General Systems 49 (2020), 574-601. https://doi.org/10.1080/03081079.2020.1786380.

\bibitem {FJ1}S. Furtado, C. R. Johnson, \textit{Efficient vectors in priority
setting methodology}, Annals of Operations Research 332 (2024), 743-764. https://doi.org/10.1007/s10479-023-05771-y.

\bibitem {FJ2}S. Furtado, C. R. Johnson, \textit{Efficiency of any weighted
geometric mean of the columns of a reciprocal matrix}, Linear Algebra and its
Applications\textit{\ }680 (2024), 83-92. https://doi.org/10.1016/j.laa.2023.10.001.

\bibitem {FJ20}S. Furtado, C. R. Johnson, \textit{Pairwise comparison matrices
with uniformly ordered efficient vectors}, International Journal of
Approximate Reasoning 173 (2024), 109265. https://doi.org/10.1016/j.ijar.2024.109265.

\bibitem {FJ3}S. Furtado, C. R. Johnson, \textit{Cycle products and efficient
vectors in reciprocal matrices}, Electronic Journal of Linear Algebra 40
(2024), 538-550. https://doi.org/10.13001/ela.2024.8543.

\bibitem {FJblock}S. Furtado, C. R. Johnson, \textit{Efficient vectors for
block perturbed consistent matrices}, SIAM Journal on Matrix Analysis and
Applications 45 (2024), 601-618. https://doi.org/10.1137/23M1580310.

\bibitem {FJ5}S. Furtado, C. R. Johnson, \textit{Efficiency analysis for the
Perron vector of a reciprocal matrix}, Applied Mathematics and Computation 480
(2024), 128913. https://doi.org/10.1016/j.amc.2024.128913.

\bibitem {FJEJOR}S. Furtado, C. R. Johnson, \textit{Efficiency analysis of
natural cardinal ranking vectors for pairwise comparisons and the universal
efficiency of the Perron geometric mean}, European Journal of Operational
Research 328 (2026), 938-945. https://doi.org/10.1016/j.ejor.2025.10.030.

\bibitem {FJ4}S. Furtado, C. R. Johnson, \textit{The complete set of efficient
vectors for a pairwise comparison matrix}, Computational and Applied
Mathematics 44 (2025), 100. https://doi.org/10.1007/s40314-024-03062-1.

\bibitem {FJ80}S. Furtado, C. R. Johnson, \textit{Best matrix interval for
consistent matrices approximating a given reciprocal matrix}, Annals of
Operations Research 358 (2026), 1199-1218. https://doi.org/10.1007/s10479-025-06898-w.

\bibitem {FJ8}S. Furtado, C. R. Johnson, \textit{Exact bounds for the
consistent matrices obtained from the efficient vectors for a reciprocal
matrix}, Linear Algebra and its Applications. https://doi.org/10.1016/j.laa.2025.12.015.

\bibitem {Gass}S. I. Gass, \textit{Tournaments, transitivity and pairwise
comparison matrices}, Journal of the Operational Research Society 49 (1998),
616-624. https://doi.org/10.1057/palgrave.jors.2600572.

\bibitem {Golany1993}B. Golany, M. Kress, \textit{A multicriteria evaluation
of methods for obtaining weights from ratio-scale matrices}, European Journal
of Operational Research 69 (1993), 210-220. https://doi.org/10.1016/0377-2217(93)90165-J.

\bibitem {Johnson}C. R. Johnson, W. B. Beine, T. J. Wang, \textit{Right-left
asymmetry in an eigenvector ranking procedure}, Journal of Mathematical
Psychology 19 (1979), 61-64. https://doi.org/10.1016/0022-2496(79)90005-1.

\bibitem {Maleki-Cs}H. Maleki, S. Zahir, \textit{A comprehensive literature
review of the rank reversal phenomenon in the Analytic Hierarchy Process},
Journal of Multi-Criteria Decision Analysis 20 (2013), 141-155. https://doi.org/10.1002/mcda.1479.

\bibitem {MAZ}J. Mazurek, R. Perzina, J. Ram\'{\i}k, D. Bartl, \textit{A
numerical comparison of the sensitivity of the geometric mean method,
Eigenvalue Method, and Best-Worst}, Method Mathematics 9 (2021), 554. https://doi.org/10.3390/math9050554.

\bibitem {Moon}J. W. Moon, \textit{Topics in Tournaments}, Holt, Rinehart and
Winston (1968).

\bibitem {saaty1977}T. L. Saaty, \textit{A scaling method for priorities in
hierarchical structures}, Journal of Mathematical Psychology 15 (1977),
234-281. https://doi.org/10.1016/0022-2496(77)90033-5.

\bibitem {Saaty-Cs}T. L. Saaty, L.G Vargas, \textit{The legitimacy of rank
reversal}, Omega 12 (1984), 513-516. https://doi.org/10.1016/0305-0483(84)90052-5.

\bibitem {TuWu}J. Tu, Z. Wu, \textit{Analytic hierarchy process rank
reversals: causes and solutions}, Annals of Operations Research 346 (2025),
1785-1809. https://doi.org/10.1007/s10479-023-05278-6.

\bibitem {Wang1-Cs}Y. M. Wang, T.M.S. Elhag, \textit{An approach to avoiding
rank reversal in AHP}, Decision Support Systems 42 (2006), 1474-1480. https://doi.org/10.1016/j.dss.2005.12.002.

\bibitem {Wang2-Cs}Y. M. Wang, Y. Luo, \textit{On rank reversal in decision
analysis}, Mathematical and Computer Modelling 49 (2009), 1221-1229. https://doi.org/10.1016/j.mcm.2008.06.019.

\bibitem {wang}H. Wang, Y. Peng, G. Kou, \textit{A two-stage ranking method to
minimize ordinal violation for pairwise comparisons,} Applied Soft Computing
106 (2021), 107287. https://doi.org/10.1016/j.asoc.2021.107287.
\end{thebibliography}
\end{document}